\documentclass[12pt,onecolumn]{IEEEtran}

\usepackage{amsmath,amssymb,psfrag,latexsym,color,cite,times,bbm,amstext,wasysym,parskip,url,graphicx,algorithm,algorithmic,tikz,float,multicol,subcaption}
\usepackage{amsthm}

\graphicspath{{figure/}{./}}

\newtheorem{theorem}{Theorem}[section]

\newtheorem{lemma}[theorem]{Lemma}
\newtheorem{definition}{Definition}[section]
\newtheorem{remark}[theorem]{Remark}

\newcommand{\calL}{\mathcal{L}}
\newcommand{\At}{\calA^{(t)}}
\newcommand{\As}{\calA^{(s)}}
\newcommand{\Ts}{T^{(s)}}
\newcommand{\Tt}{T^{(t)}}
\newcommand{\calW}{\mathcal{W}}
\newcommand{\calG}{\mathcal{G}}
\newcommand{\calV}{\mathcal{V}}
\newcommand{\calE}{\mathcal{E}}
\newcommand{\calP}{\mathcal{P}}

\newcommand{\calA}{\mathcal{A}}
\newcommand{\xp}{x^{(p)}}
\newcommand{\yst}{y^{(s,t)}}

\DeclareMathOperator*{\argmin}{arg\,min}

\makeatletter
\newcommand\fs@norules{\def\@fs@cfont{\bfseries}\let\@fs@capt\floatc@ruled
  \def\@fs@pre{}%
  \def\@fs@post{}%
  \def\@fs@mid{\kern3pt}%
  \let\@fs@iftopcapt\iftrue}
\makeatother
\floatstyle{norules}
\restylefloat{algorithm}

\def\spacingset#1{\def\baselinestretch{#1}\small\normalsize}
\spacingset{1}

\title{Lasso formulation of the shortest path problem}
\author{Anqi Dong, Amirhossein Taghvaei, Tryphon T.\ Georgiou
\thanks{A. Dong, A. Taghvaei, and T. T. Georgiou are with the Department of Mechanical and Aerospace Engineering, University of California, Irvine, CA 92697-3975; \{anqid2,ataghvae,tryphon\}@uci.edu}}

\begin{document}
	\maketitle
	
	\begin{abstract}
		The  shortest path problem is formulated as an $l_1$-regularized regression problem, known as lasso. Based on this formulation, a connection is established between Dijkstra's shortest path algorithm and the least angle regression (LARS) for the lasso problem.  Specifically, the solution path of the lasso problem, obtained by varying the regularization parameter from infinity to zero (the regularization path), corresponds to shortest path trees that appear in the bi-directional Dijkstra algorithm. Although Dijkstra's algorithm and the LARS formulation provide exact solutions, they become impractical  when the size of the graph is exceedingly large. To overcome this issue, the alternating direction method of multipliers (ADMM) is proposed to solve the lasso formulation. The resulting algorithm produces good and fast approximations of the shortest path by sacrificing exactness that may not be absolutely essential in many applications. Numerical experiments are provided to illustrate the performance of the proposed approach. 
	\end{abstract}

\section{Introduction}

The problem of finding the shortest path between two vertices in a graph has a long history~\cite{wiener1873ueber},\cite{tarry1895probleme} with a wide range of applications Waxman~\cite{waxman1988routing} Mortensen et al.~\cite{mortensen1995intelligent}. The classical algorithm to determine a shortest path is due to  
Edsger W. Dijkstra \cite{dijkstra1959note}. Since Dijkstra's early work, a variety of alternative methods to identify a shortest path have been developed to reduce  complexity and improve speed~\cite{bast2016route,ahuja1990faster,van1976design,fredman1987fibonacci}. However, it is often the case that finding a shortest path is not absolutely essential, especially in graphs of considerably large sizes, while a reasonably short path may suffice~\cite{waxman1988routing,potamias2009fast}. Motivated by such considerations and inspired by success of convex optimization to address large-scale problems~\cite{boyd2004convex,boyd2011distributed}, we introduce a formulation of the shortest path problem as an $l_1$-regularized regression, known as the ``lasso'' (Least Absolute Shrinkage and Selection Operator) problem~\cite{tibshirani1996regression}. 

Specifically, in this paper, we discuss two novel and important implications of the lasso formulation which constitute our main contributions. 

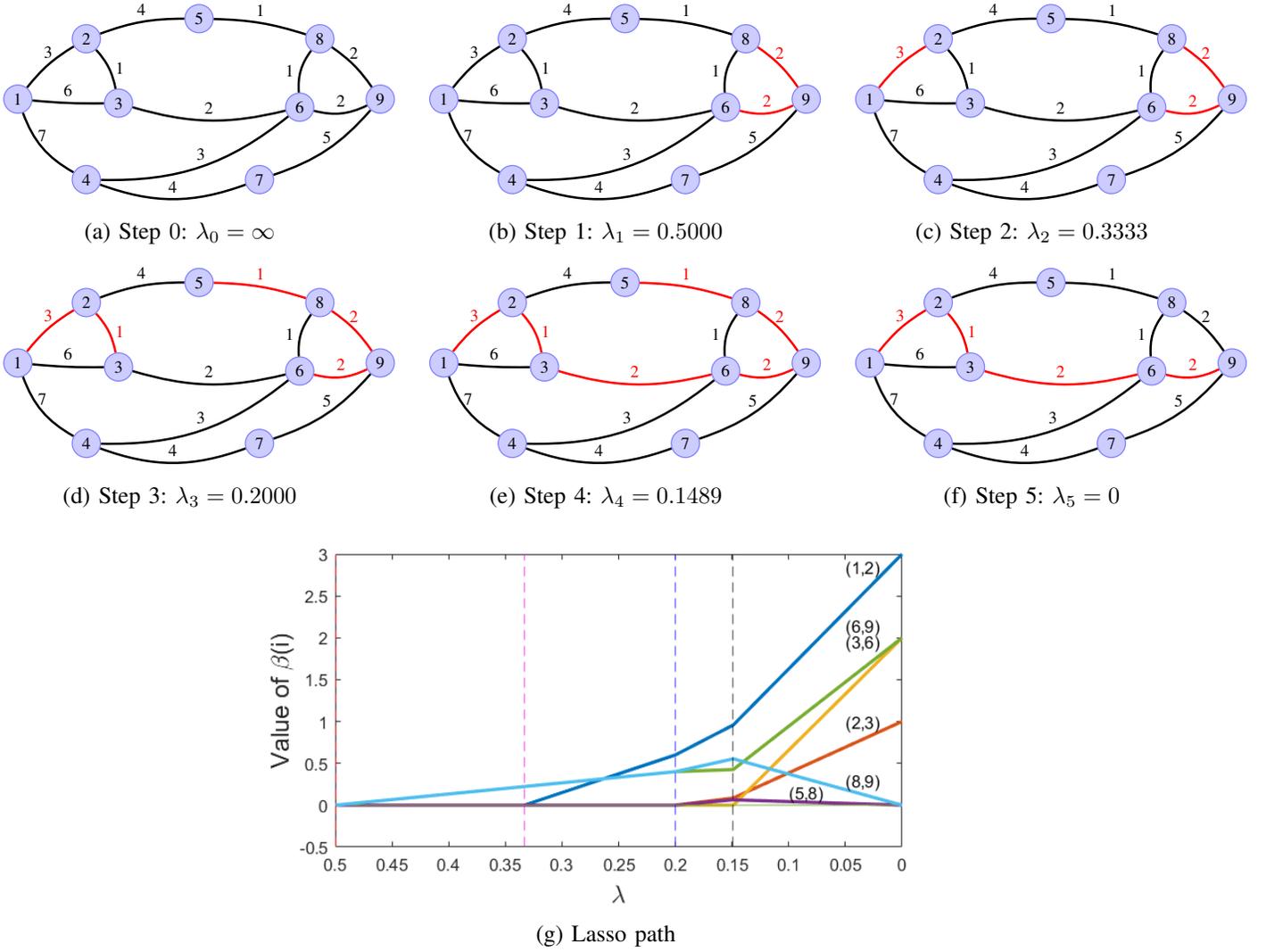
\begin{figure}[H]
	\begin{subfigure}{0.3\columnwidth}
		\centering
		\scalebox{0.6}{
			\begin{tikzpicture}
			[process/.style={circle,draw=blue!50,fill=blue!20,thick,
				inner sep=0pt,minimum size=7mm},
			target/.style={rectangle,draw=black!50,fill=black!20,thick,
				inner sep=0pt,minimum size=8mm}]
			\node (n1) at (0,0) [process] {1};
			\node (n2) at (1.7,1.5) [process] {2};
			\node (n3) at (2.5,-0.1) [process] {3};
			\node (n4) at (1.7,-2) [process] {4};
			\node (n5) at (4.5,2) [process] {5};
			\node (n6) at (7,-0.2) [process] {6};
			\node (n7) at (6,-2) [process] {7};
			\node (n8) at (7.5,1.5) [process] {8};
			\node (n9) at (9,0) [process] {9};  
			
			\draw[ultra thick] (n1) to [bend left=10] node[midway,above](){3} (n2);
			\draw[ultra thick] (n1) to [bend right=3] node[midway,above](){6} (n3);
			\draw[ultra thick] (n1) to [bend right=20] node[midway,above](){7} (n4);
			\draw[ultra thick] (n2) to [bend left=10] node[midway,above](){4} (n5);
			\draw[ultra thick] (n2) to [bend left=15] node[midway,right](){1} (n3);
			\draw[ultra thick] (n3) to [bend right=15] node[midway,above](){2} (n6);
			\draw[ultra thick] (n4) to [bend right=20] node[midway,above](){3} (n6);
			\draw[ultra thick] (n4) to [bend right=20] node[midway,above](){4} (n7);
			\draw[ultra thick] (n5) to [bend left=8] node[midway,above](){1} (n8);
			\draw[ultra thick] (n6) to [bend left=20] node[midway,left](){1} (n8);
			\draw[ultra thick] (n6) to [bend right=20] node[midway,above](){2} (n9);
			\draw[ultra thick] (n7) to [bend right=15] node[midway,above](){5} (n9);
			\draw[ultra thick] (n8) to [bend left=10] node[midway,above](){2} (n9);
			\end{tikzpicture}}
		\caption{Step $0$: $\lambda_{0} = \infty$} \label{fig:a}
	\end{subfigure}\hspace*{\fill}
	\begin{subfigure}{0.3\columnwidth}
		\centering
		\scalebox{0.6}{
			\begin{tikzpicture}
			[process/.style={circle,draw=blue!50,fill=blue!20,thick,
				inner sep=0pt,minimum size=7mm},
			target/.style={rectangle,draw=black!50,fill=black!20,thick,
				inner sep=0pt,minimum size=8mm}]
			\node (n1) at (0,0) [process] {1};
			\node (n2) at (1.7,1.5) [process] {2};
			\node (n3) at (2.5,-0.1) [process] {3};
			\node (n4) at (1.7,-2) [process] {4};
			\node (n5) at (4.5,2) [process] {5};
			\node (n6) at (7,-0.2) [process] {6};
			\node (n7) at (6,-2) [process] {7};
			\node (n8) at (7.5,1.5) [process] {8};
			\node (n9) at (9,0) [process] {9};  
			
			\draw[ultra thick] (n1) to [bend left=10] node[midway,above](){3} (n2);
			\draw[ultra thick] (n1) to [bend right=3] node[midway,above](){6} (n3);
			\draw[ultra thick] (n1) to [bend right=20] node[midway,above](){7} (n4);
			\draw[ultra thick] (n2) to [bend left=10] node[midway,above](){4} (n5);
			\draw[ultra thick] (n2) to [bend left=15] node[midway,right](){1} (n3);
			\draw[ultra thick] (n3) to [bend right=15] node[midway,above](){2} (n6);
			\draw[ultra thick] (n4) to [bend right=20] node[midway,above](){3} (n6);
			\draw[ultra thick] (n4) to [bend right=20] node[midway,above](){4} (n7);
			\draw[ultra thick] (n5) to [bend left=8] node[midway,above](){1} (n8);
			\draw[ultra thick] (n6) to [bend left=20] node[midway,left](){1} (n8);
			\draw[ultra thick][red] (n6) to [bend right=20] node[midway,above](){2} (n9);
			\draw[ultra thick] (n7) to [bend right=15] node[midway,above](){5} (n9);
			\draw[ultra thick][red] (n8) to [bend left=10] node[midway,above](){2} (n9);
			\end{tikzpicture}}
		\caption{Step 1: $\lambda_{1} = 0.5000$} \label{fig:b}
	\end{subfigure}\hspace*{\fill}
	\begin{subfigure}{0.3\columnwidth}
		\centering
		\scalebox{0.6}{
			\begin{tikzpicture}
			[process/.style={circle,draw=blue!50,fill=blue!20,thick,
				inner sep=0pt,minimum size=7mm},
			target/.style={rectangle,draw=black!50,fill=black!20,thick,
				inner sep=0pt,minimum size=8mm}]
			\node (n1) at (0,0) [process] {1};
			\node (n2) at (1.7,1.5) [process] {2};
			\node (n3) at (2.5,-0.1) [process] {3};
			\node (n4) at (1.7,-2) [process] {4};
			\node (n5) at (4.5,2) [process] {5};
			\node (n6) at (7,-0.2) [process] {6};
			\node (n7) at (6,-2) [process] {7};
			\node (n8) at (7.5,1.5) [process] {8};
			\node (n9) at (9,0) [process] {9};  
			
			\draw[ultra thick][red] (n1) to [bend left=10] node[midway,above](){3} (n2);
			\draw[ultra thick] (n1) to [bend right=3] node[midway,above](){6} (n3);
			\draw[ultra thick] (n1) to [bend right=20] node[midway,above](){7} (n4);
			\draw[ultra thick] (n2) to [bend left=10] node[midway,above](){4} (n5);
			\draw[ultra thick] (n2) to [bend left=15] node[midway,right](){1} (n3);
			\draw[ultra thick] (n3) to [bend right=15] node[midway,above](){2} (n6);
			\draw[ultra thick] (n4) to [bend right=20] node[midway,above](){3} (n6);
			\draw[ultra thick] (n4) to [bend right=20] node[midway,above](){4} (n7);
			\draw[ultra thick] (n5) to [bend left=8] node[midway,above](){1} (n8);
			\draw[ultra thick] (n6) to [bend left=20] node[midway,left](){1} (n8);
			\draw[ultra thick][red] (n6) to [bend right=20] node[midway,above](){2} (n9);
			\draw[ultra thick] (n7) to [bend right=15] node[midway,above](){5} (n9);
			\draw[ultra thick][red] (n8) to [bend left=10] node[midway,above](){2} (n9);
			\end{tikzpicture}}
		\caption{Step 2: $\lambda_{2} = 0.3333$} \label{fig:c}
	\end{subfigure}

    \medskip
	\begin{subfigure}{0.3\columnwidth}
		\centering
		\scalebox{0.6}{
			\begin{tikzpicture}
			[process/.style={circle,draw=blue!50,fill=blue!20,thick,
				inner sep=0pt,minimum size=7mm},
			target/.style={rectangle,draw=black!50,fill=black!20,thick,
				inner sep=0pt,minimum size=8mm}]
			\node (n1) at (0,0) [process] {1};
			\node (n2) at (1.7,1.5) [process] {2};
			\node (n3) at (2.5,-0.1) [process] {3};
			\node (n4) at (1.7,-2) [process] {4};
			\node (n5) at (4.5,2) [process] {5};
			\node (n6) at (7,-0.2) [process] {6};
			\node (n7) at (6,-2) [process] {7};
			\node (n8) at (7.5,1.5) [process] {8};
			\node (n9) at (9,0) [process] {9};  
			
			\draw[ultra thick][red] (n1) to [bend left=10] node[midway,above](){3} (n2);
			\draw[ultra thick] (n1) to [bend right=3] node[midway,above](){6} (n3);
			\draw[ultra thick] (n1) to [bend right=20] node[midway,above](){7} (n4);
			\draw[ultra thick] (n2) to [bend left=10] node[midway,above](){4} (n5);
			\draw[ultra thick][red] (n2) to [bend left=15] node[midway,right](){1} (n3);
			\draw[ultra thick] (n3) to [bend right=15] node[midway,above](){2} (n6);
			\draw[ultra thick] (n4) to [bend right=20] node[midway,above](){3} (n6);
			\draw[ultra thick] (n4) to [bend right=20] node[midway,above](){4} (n7);
			\draw[ultra thick][red] (n5) to [bend left=8] node[midway,above](){1} (n8);
			\draw[ultra thick] (n6) to [bend left=20] node[midway,left](){1} (n8);
			\draw[ultra thick][red] (n6) to [bend right=20] node[midway,above](){2} (n9);
			\draw[ultra thick] (n7) to [bend right=15] node[midway,above](){5} (n9);
			\draw[ultra thick][red] (n8) to [bend left=10] node[midway,above](){2} (n9);
			\end{tikzpicture}}
		\caption{Step 3: $\lambda_{3} = 0.2000$} \label{fig:d}
	\end{subfigure}\hspace*{\fill}
	\begin{subfigure}{0.3\columnwidth}
		\centering
		\scalebox{0.6}{
			\begin{tikzpicture}
			[process/.style={circle,draw=blue!50,fill=blue!20,thick,
				inner sep=0pt,minimum size=7mm},
			target/.style={rectangle,draw=black!50,fill=black!20,thick,
				inner sep=0pt,minimum size=8mm}]
			\node (n1) at (0,0) [process] {1};
			\node (n2) at (1.7,1.5) [process] {2};
			\node (n3) at (2.5,-0.1) [process] {3};
			\node (n4) at (1.7,-2) [process] {4};
			\node (n5) at (4.5,2) [process] {5};
			\node (n6) at (7,-0.2) [process] {6};
			\node (n7) at (6,-2) [process] {7};
			\node (n8) at (7.5,1.5) [process] {8};
			\node (n9) at (9,0) [process] {9};  
			
			\draw[ultra thick][red] (n1) to [bend left=10] node[midway,above](){3} (n2);
			\draw[ultra thick] (n1) to [bend right=3] node[midway,above](){6} (n3);
			\draw[ultra thick] (n1) to [bend right=20] node[midway,above](){7} (n4);
			\draw[ultra thick] (n2) to [bend left=10] node[midway,above](){4} (n5);
			\draw[ultra thick][red] (n2) to [bend left=15] node[midway,right](){1} (n3);
			\draw[ultra thick][red] (n3) to [bend right=15] node[midway,above](){2} (n6);
			\draw[ultra thick] (n4) to [bend right=20] node[midway,above](){3} (n6);
			\draw[ultra thick] (n4) to [bend right=20] node[midway,above](){4} (n7);
			\draw[ultra thick][red] (n5) to [bend left=8] node[midway,above](){1} (n8);
			\draw[ultra thick] (n6) to [bend left=20] node[midway,left](){1} (n8);
			\draw[ultra thick][red] (n6) to [bend right=20] node[midway,above](){2} (n9);
			\draw[ultra thick] (n7) to [bend right=15] node[midway,above](){5} (n9);
			\draw[ultra thick][red] (n8) to [bend left=10] node[midway,above](){2} (n9);
			\end{tikzpicture}}
		\caption{Step 4: $\lambda_{4} = 0.1489$} \label{fig:e}
	\end{subfigure}\hspace*{\fill}
	\begin{subfigure}{0.3\columnwidth}
		\centering
		\scalebox{0.6}{
			\begin{tikzpicture}
			[process/.style={circle,draw=blue!50,fill=blue!20,thick,
				inner sep=0pt,minimum size=7mm},
			target/.style={rectangle,draw=black!50,fill=black!20,thick,
				inner sep=0pt,minimum size=8mm}]
			\node (n1) at (0,0) [process] {1};
			\node (n2) at (1.7,1.5) [process] {2};
			\node (n3) at (2.5,-0.1) [process] {3};
			\node (n4) at (1.7,-2) [process] {4};
			\node (n5) at (4.5,2) [process] {5};
			\node (n6) at (7,-0.2) [process] {6};
			\node (n7) at (6,-2) [process] {7};
			\node (n8) at (7.5,1.5) [process] {8};
			\node (n9) at (9,0) [process] {9};  
			
			\draw[ultra thick][red] (n1) to [bend left=10] node[midway,above](){3} (n2);
			\draw[ultra thick] (n1) to [bend right=3] node[midway,above](){6} (n3);
			\draw[ultra thick] (n1) to [bend right=20] node[midway,above](){7} (n4);
			\draw[ultra thick] (n2) to [bend left=10] node[midway,above](){4} (n5);
			\draw[ultra thick][red] (n2) to [bend left=15] node[midway,right](){1} (n3);
			\draw[ultra thick][red] (n3) to [bend right=15] node[midway,above](){2} (n6);
			\draw[ultra thick] (n4) to [bend right=20] node[midway,above](){3} (n6);
			\draw[ultra thick] (n4) to [bend right=20] node[midway,above](){4} (n7);
			\draw[ultra thick] (n5) to [bend left=8] node[midway,above](){1} (n8);
			\draw[ultra thick] (n6) to [bend left=20] node[midway,left](){1} (n8);
			\draw[ultra thick][red] (n6) to [bend right=20] node[midway,above](){2} (n9);
			\draw[ultra thick] (n7) to [bend right=15] node[midway,above](){5} (n9);
			\draw[ultra thick] (n8) to [bend left=10] node[midway,above](){2} (n9);
			\end{tikzpicture}}
		\caption{Step 5: $\lambda_{5} = 0$} \label{fig:f}
	\end{subfigure}
	
   \medskip
   \begin{subfigure}{1\columnwidth}
	\centering
	\includegraphics[width=0.6\columnwidth]{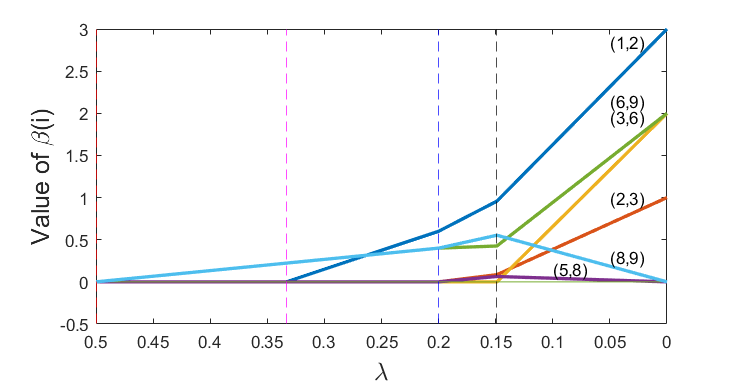}
	\caption{Lasso path} \label{fig:lassopath}
   \end{subfigure}
	\caption{The shortest path problem on the Nicholson's graph between vertices 1 and 9. The figures depict the solution path of the lasso formulation of the shortest path problem~\eqref{eq:lasso}, obtained via the LARS algorithm~\ref{alg:LARS}. Subfigures (a-b-c-d-e-f) highlight the steps of the LARS algorithm corresponding to different values of $\lambda$; the edges that belong to the active set are highlighted with red. Subfigure (g) displays the value of each component of the incidence vector $\beta$ as a function of $\lambda$. } 	\label{fig:example}
\end{figure}

(i) We provide a rather surprising connection between Dijkstra's algorithm  and the solution path of the lasso problem; we show that the solution path of the lasso problem generates shortest path trees that appear in Dijkstra's algorithm. The connection is interesting as the lasso solution path is based on analytical arguments, invoking KKT conditions, unlike the Dijkstra's algorithm that is akin to Dynamic Programming, cf.\ Figure~\ref{fig:example}.

(ii) On the practical side, we consider the shortest path problem on graphs with large size  and propose to utilize the ADMM method to obtain approximate shortest path solutions. Moreover, the ADMM method can be implemented in a distributed manner, and has the flexibility to be initialized with a rough approximation of the shortest path (if one such is available) for faster convergence; this option arise in cases where a graph undergoes slight variation from an earlier one where a short path is available, cf.\ Figure~\ref{fig:Convergence-lasso-RG-ADMM} and Figure ~\ref{fig:ConvergenceADMM-pika-ADMM}

 Here is the outline of the paper: The necessary preliminary definitions and notations are introduced in Section~\ref{sec:prelim}. The Lasso formulation is presented in Section~\ref{sec:Lasso}. The analysis of solution path of the lasso appears in Section~\ref{sec:LARS}. The  connection with Dijkstra's algorithm appears in Section~\ref{sec:relation}. The ADMM algorithm and its application to shortest path problem to large graphs is presented in Section~\ref{sec:proximal}.
 
%In the later sections of this paper, we will first formulate the well-defined shortest path problem into a $l_{1}$ regularization problem with the incidence matrix of the graph $\mathcal{G}$ and solve this problem by the LARS algorithm. Then by our derivation of the $t^{join}$ and $t^{cross}$ in LARS, we prove that the LARS algorithm is equal to Dijkstra algorithm. Finally, we apply ADMM to solve the shortest path problem in large scale graphs.

\section{Preliminaries } \label{sec:prelim}

\subsection{Graph theoretic notations and definitions}
Consider a weighted undirected  graph $\mathcal{G}=(\calV,\calE,\calW)$ with vertex set $\mathcal{V} = \{1, \dots , n\}$ and a set of edges $\mathcal{E} = \{e_{1}, \dots, e_{m}\}$ with corresponding positive weights in $\calW=\{w_1,\ldots,w_m\}$. The graph is assumed to be connected and simple (i.e., with no self-loops or multi-edges). The edge that connects the nodes $i$ and $j$ is also denoted by the pair of nodes as $(i,j)$. Although the graph is undirected, we assign an arbitrary but fixed orientation to each edge. 
%The cardinality of 
%the vertex set and the edge set is denoted as $|\mathcal{V}|=n$ and $|\mathcal{E}|=m$.
%Vertices $v_{i}$ and $v_{j}$ are adjacent when 
%\begin{align*}
%    e_{ij}=\{v_{i},\ v_{j}\}, \ e_{ij} \in \mathcal{E}
%\end{align*}

The incidence matrix of the graph, denoted by $D$, is a $n\times m$ matrix whose $(i,j)$-th entry is defined according to
\begin{equation*}
[D]_{ij} =
\begin{cases}
+1 \ \ \  &\text{if vertex $i$ is at the tail of edge $e_{j}$},\\
-1  \ \ \  &\text{if vertex $i$ is at the head of edge $e_{j}$},\\
0  \ \ \  &\text{otherwise}.
\end{cases}  
\end{equation*}
The graph Laplacian matrix, denoted by $L$, is a $n\times n$ matrix defined as
\begin{align*}
L = DWD^T 
\end{align*}
where $W=\text{diag}(w_1,\ldots,w_m)$ is the diagonal matrix formed by the weights. 

A path from vertex $s$  to vertex $t$ is a sequence of connected edges $p=\{(v_0,v_1),(v_1,v_2),\ldots, (v_{l-1},v_l)\} \subset \calE$ that originates  
at $v_0=s$ and ends at $v_l=t$. The incidence vector $\xp$ of a path $p$ is $m$-dimensional vector defined as follows. The $i$-th entry of $\xp$ is zero, if the path does not contain the edge $e_i$. If the path contains $e_i$, then the corresponding entry is either $+1$ or $-1$. The sign is positive/negative depending on whether the direction of the path agrees/disagrees with direction of the edge. The length of the path is the sum of the weights of the edges that belong to the path, i.e.,
\begin{equation*}
\text{length}(p) \triangleq \sum_{e_i \in p} w_i = \|W\xp\|_1
\end{equation*} 
where  $\|\cdot\|_1$ denotes the $l_1$-norm.

%\section{Preliminaries}
%\subsection{Graph theoretic notations and definitions:}
%The graph $\mathcal{G}$ is defined as $\mathcal{G} = (\mathcal{V},\; \mathcal{E})$. $\mathcal{V} = \{v_{1} \dots v_{n}\}$ is the vertices set which contains all the $n$ vertices, $\mathcal{E} = \{e_{1} \dots e_{m}\}$ is the edge set. The cardinality of 
%the vertex set and the edge set is denoted as $|\mathcal{V}|=n$ and $|\mathcal{E}|=m$.
%Vertices $v_{i}$ and $v_{j}$ are adjacent when 
%\begin{align*}
%e_{ij}=\{v_{i},\ v_{j}\}, \ e_{ij} \in \mathcal{E}
%\end{align*}
%The graph $\mathcal{G}$ is connected if for every pair of vertices in $\mathcal{V}$ there exists a path between them.
%The incidence matrix and the graph Laplacian is defined as following:\\
%
%\begin{definition}
%	The incidence matrix is a $n\times m$ matrix with $n$ rows and $m$ columns obeying the  following rules
%	\begin{align*}
%	D = [d_{ij}], \ \ \ d_{ij} =
%	\begin{cases}
%	-1 \ \ \  &\textit{if $v_{i}$ is the tail of $e_{j}$},\\
%	1  \ \ \  &\textit{if $v_{i}$ is the head of $e_{j}$},\\
%	0  \ \ \  &\textit{otherwise}.
%	\end{cases}.    
%	\end{align*}
%\end{definition}
%\vspace{0.1in}
%
%\begin{definition}
%	The weights matrix $W$ is a diagonal matrix such that the diagonal element $w_{ii}$ is the length of the edge $e_{i}$. 
%\end{definition}

\begin{definition}\label{def:P}
	When $\calG$ is a tree with root at vertex $1$, the path matrix $P$ is  defined to be the $(n-1)\times (n-1)$ matrix whose $i$-th column is the incidence vector of the path from vertex $i$ to the root $1$, for $i=2,3,\ldots,n$.   
%	are equal to the incidence The rows of the path matrix correspond to the $(n-1)$ edges of the tree and the columns of the matrix correspond to the $(n-1)$ vertices except the root. The $(i,j)$-th entry of $P$ is equal to zero if the edge $e_i$ belongs to the path of vertex $j$ to the root the existence of edge $e_{i}$ in the path of $v_{j}$ to the root of $\mathcal{G}$ with it's direction.
\end{definition}

The pseudo-inverse of the incidence matrix of a tree admits a  simple representation in terms of the path matrix~\cite[Theorem 2.10, Lemma 2.15]{bapat2010graphs}. 
\begin{lemma}\label{lem:incidence}
	Let $D$ be the incidence matrix of a tree with vertices $\calV=\{1,2,\ldots,n\}$ and root vertex $1$. Then, the pseudo-inverse of $D$ is given by
	\begin{equation}
	D_{}^{+} =
	\begin{bmatrix}
	-\frac{1}{n}P \mathbbm{1} & PJ
	\end{bmatrix}_{(n-1)\times n}
	\end{equation}
	where 
	\[
	\mathbbm{1} := [1,1,\ldots,1]^T_{n-1}
	\]
	 is $(n-1)$-dimensional (column) vector of ones, 
	 \[J=(I-\frac{1}{n} \mathbbm{1} \mathbbm{1}^T)\]
	 is an orthogonal projection with null space spanned by $\mathbbm{1}$,
	   $I$ is identity matrix of size $(n-1)\times (n-1)$, and $P$ is the $(n-1)\times(n-1)$ path matrix of the tree, defined in~\ref{def:P}.  
\end{lemma}

%\begin{definition}
%	Given a connected, undirected graph $\mathcal{G}$, a shortest path tree (\textit{SPT}) rooted at vertex $v_{s}$ is a spanning tree $\mathcal{T}$ of $\mathcal{G}$, such that the path distance from root $v_{s}$ to any other vertex $v_{i}$  in $\mathcal{T}$ is the shortest path distance from $v_{s}$ to $v_{i}$ in $\mathcal{G}$.
%\end{definition}

\subsection{Shortest path problem}
%We will bring up the shortest path problem and the lasso problem in this section. By connecting the path set and the active set in this two problems, we are able to show that the shortest path problem can be also solved by lasso. Moreover, the new algorithm we propose in the later chapter shows that after some significant simplify, this specific type of lasso is equal to  the Bidirectional Dijkstra.\\

%\subsection{The shortest path problem}
%A path between two vertices $s$ and $t$ is a sequence of edges $\{(v_0,v_1),(v_1,v_2),\ldots, (v_{l-1},v_l)\} \subset \calE$ that originates  
%at $v_0=s$ and ends at $v_l=t$. 
Let $\mathcal{P}_{s,t}$ denote the set of all paths between  $s$ and $t$. This set is non-empty because the graph is connected. 
The {\em shortest path problem} is to find a path between  $s$ and $t$ with minimum length, mathematically formulated as finding
\begin{equation}\label{eq:shortest-path-problem}
{\rm arg}\min_{p\in \mathcal P_{s,t}}~\text{length}(p).
\end{equation}
The minimum value is the distance between $s$ and $t$, and the minimizing path is the shortest path between $s$ and $t$.  
%The length of a path $p\in \mathcal P_{s,t}$ is the sum of the weights of the edges in the path, i.e. 
%\begin{equation}
%\text{length}(p) = \sum_{e_i \in p} w_i
%\end{equation} 

%\subsection{The shortest path problem}
%The path $\mathcal{P}$ of two vertices $v_{s},\ v_{t} \in \mathcal{V}$ is a sequence of edges \cite{pohl1969bidirectional}, $v_{s}$ is the starting vertex and $v_{t}$ is the targeting vertex so that we can define the $\mathcal{P}$ as below
%\begin{align*}
%\mathcal{P} &= \{e_{1},\ e_{2},\ \dots,\ e_{k}\},\\ 
% \text{s.t.}\ \ &e_{i} \cap e_{i+1} \neq \emptyset,\\
% & e_{i} \neq e_{i+1},\\
% & v_{s} \in e_{1}, \ v_{t} \in e_{k}\\
% & \forall \ 1\leq i \leq k-1.
%\end{align*}
%
%Therefore, we can formulate the shortest path problem for weighted graph.
%For the weighted graph $\mathcal{G} = \{\mathcal{V},\ \mathcal{E}\}$, since each edge $e_{i}$ has a unique weight $w(e_{i})\geq 0$, the total length of the path can be written as. 
%\begin{align}
%w(\mathcal{P}) = \sum_{1}^{k} w(e_{i}), \ e_{i} \in \mathcal{P}
%\end{align}
%
%For all paths between $v_{s}$ and $v_{t}$, the shortest path $\hat{\mathcal{P}}$ should be the path with minimum weight $w(\mathcal{P})$.

\subsection{Dijkstra's algorithm}\label{sec:Dikstra}
Dijkstra's algorithm to find the shortest path from  $v_s$ to $v_t$ involves the following variables:
\begin{itemize}
	\item \textit{dist}: an array of distances from the root vertex $v_s$ to all the other vertices in the graph. 
	\item $\mathcal S$: the set of visited vertices.
	\item $\mathcal Q$: the queue of vertices to be visited. 
\end{itemize}
The algorithm begins with initial $\infty$ values for the distances   and improve the distance step by step as follows
%\AHT{can you write the algorithm is pseudo-code? (Done)}
\begin{algorithm}[H]
	\caption{\bf Dijkstra's algorithm}
	\begin{algorithmic}[1]
		\renewcommand{\algorithmicrequire}{\textbf{Input:}}
		\renewcommand{\algorithmicensure}{\textbf{Output:}}
		\REQUIRE source vertex $v_s$ and target  vertex $v_t$. %\hspace{1.4in}// \textit{dist}[$v_{i}$] is the distance from $v_{i}$ to $v_{s}. $
		\ENSURE  the shortest path and the length of the path.
		\STATE   $dist^{(0)}$[$v_s$]$=0$, $dist^{(0)}$[$v_{i}$]$=\infty, \forall v_i \neq v_s$, $\mathcal S = \emptyset$, $\mathcal Q=\calV$% \hspace{1.1in} // $S$ is the set of visited vertices.
		\WHILE{$v_t \notin S$} 	
		\STATE pick $u$ from $Q$ with minimum distance: %\hspace{1.1in}// $Q$ is the queue which initially equal to $\mathcal{V}$.\\ 
               \hspace{0.3in}$u = \argmin\limits_{v \in Q} dist[v]$
		\STATE Remove $u$ from $Q$: $Q \leftarrow Q\setminus \{u\}$
		\STATE Add $u$ to $S$: $S \leftarrow S \cup \{u\}$
		\FOR { $ v_{i} \in \textit{neighbors}[u]$}
		      \IF {\textit{dist}[$v_{i}$] $>$ \textit{dist}[$u$]$+w_{u,v_{i}}$} 
		    %   \hspace{1.5in} // $w_{u,v_{i}}$ is the weight on edge $(u,v_{i})$.\\
		      % \hspace{0.4in}then\\
               \STATE $\textit{dist}[v_{i}] \leftarrow \textit{dist}[u]+w_{u,v_{i}}$
               \ENDIF 
               \ENDFOR
		\ENDWHILE
		\RETURN \textit{dist}[$v_{t}$]
	\end{algorithmic}
	\label{alg:Dijkstra}
\end{algorithm}
%where \textit{dist} is an array of distance from the root vertex $s$ to all the other vertices in the graph $\calG$, $S$ is the set of visited vertices and $Q$ is the queue which initially contains all vertices.\\[0.1in] 
%In each step of the Dijkstra's algorithm, the  \textit{dist} array  is updated to reach its correct value at the end. Each step is as follows: 
%\begin{itemize}
%	\item While the target vertex $t$ is not in the visited vertices set $S$, select the element $v_{i}$ of $Q$ with the minimum distance.
%	\item Add $v_{i}$ to $S$ and marked it as visited.
%	\item For all the $v_{j}$ that are adjacent to $v_{i}$, if \textit{dist}[$v_{j}$] is larger than \textit{dist}[$v_{i}$]$+w(v_{i},v_{j})$, then
%	\begin{align*}
%	\textit{dist}[v_{j}] \leftarrow \textit{dist}[v_{i}]+w(v_{i},v_{j})
%	\end{align*} 
%\end{itemize}
The essential feature of the Dijkstra’s algorithm is that, the algorithm iteratively constructs the shortest path tree that is rooted at $s$, to all the visited vertices before reaching the target $t$. Later in Section~\ref{sec:relation}, we show that such a feature is also observed in the Lasso formulation of the shortest path problem.  
%The Dijkstra's algorithm creates a tree of shortest paths from $s$ to all the visited vertices in $S$, which is defined as the shortest path tree. So when the $t$ is added to the $S$, the algorithm can ensure that the \textit{dist}[$t$] is already the length of the shortest path. In another words, adding $v_{i}$ to $S$ means that the shortest path from $s$ to $v_{i}$ has already be found.

%\newpage
\section{Lasso formulation of the shortest path problem}\label{sec:Lasso}

%\subsection{Incidence vector and incidence matrix} 
%Each path in the graph is identified with a $m$-dimensional vector called incidence vector. The $i$-th entry of the vector is equal to $+1$ or $-1$ if the edge $e_i$ belongs to the path. Otherwise the entry is zero. The sign for the nonzero entry is positive if the path follow the direction of the edge, and negative otherwise. The definition of the incidence vector $x \in \mathbb R^m$ for a path $p$ is as follows:
%\begin{align}\label{betadef}
%x_i =  
%\begin{cases}
%\pm 1, & \text{if} \ \ e_{i} \in \mathcal{P} \\
%0, & \text{if} \ \ e_{i} \notin \mathcal{P}
%\end{cases}
%\end{align}  
\subsection{Linear programming formulation}
The shortest path problem~\eqref{eq:shortest-path-problem} can be formulated as a linear programming problem on the incidence vectors. 
The condition that a path, or in general a collection of edges, forms a path between $s$ and $t$ can be expressed as a linear constraint on the incidence vector: 
\begin{equation}\label{eq:Dxy}
p \in \calP_{s,t} \quad \Rightarrow\quad D\xp  = \yst
\end{equation}
where $\yst \in \mathbb R^n$ is defined according to 
\begin{align*}
\yst(i) = 
\begin{cases}
+1, & \text{if}\ i = s, \\
-1,  & \text{if}\ i =t, \\ 
0, &\text{otherwise}.
\end{cases} \ \ \ 1 \leq i \leq n
\end{align*}
The justification for~\eqref{eq:Dxy} is straightforward. 
%The $i$-th element of $D\xp$ is $(D\xp)_i=\sum_{j \in \calN(i)}D_{i,j} \xp_j$, where $\calN(i)$ is the set of all neighbours of vertex $i$. The entry is clearly zero when the path does not intersect vertex $i$. Otherwise, there are precisely two edges of the path that intersect $i$ when $i\neq s,t$ and one edge of $i=s,t$.    When $i \neq s,t$, the path crosses $i$ or   does not belong to the path.  
For any two connecting edges $e_i=(v_1,v_2)$ and $e_j=(v_2,v_3)$, the summation of the $i$-th and $j$-th columns of $D$ is equal to $y^{(v_1,v_3)}$, which corresponds to the edge $(v_1,v_3)$. Likewise, an additional column corresponding to the edge $(v_3,v_4)$, yields $y^{(v_1,v_4)}$.  Therefore, $D \xp$, results in the summation of all columns corresponding to a set of connecting edges $p=\{(s,v_1),(v_1,v_2),\ldots, (v_{l-1},t)\}$, and this is $y^{(s,t)}$.\\[-.05in]

\begin{remark}An alternative justification can be provided by noting that closed cycles, i.e., paths that begin and end at the same node, form a basis for the null space of the incidence matrix of the graph \cite{bapat2010graphs,mesbahi2010graph}. If we attach a {\em virtual direct link (i.e., a new edge)} between vertices $s$ and $t$, we need to update the incidence matrix to
$\left[\begin{matrix} D& -y^{(s,t)}\end{matrix}\right]$ so that this virtual edge is included. Now a path from $s$ to $t$ ``closes'' into a cycle by including this extra virtual edge. Any cycle that includes the virtual edge corresponds to a null vector of $\left[\begin{matrix} D& -y^{(s,t)}\end{matrix}\right]$
with a $1$ as the last entry (indicating that the virtual edge is included), and therefore, to a solution of
\[
\left[\begin{matrix} D& -y^{(s,t)}\end{matrix}\right] \left[\begin{matrix} x^{(p)}\\ 1\end{matrix}\right]=0.
\]
This is precisely \eqref{eq:Dxy}, while the first component $x^{(p)}$ of the solution vector corresponds to a sought path from $s$ to $t$.
%The linear constraint~\eqref{eq:Dxy} can be also understood in terms of cycles and kernel of incidence matrix. Consider a new graph where the edge $(t,s)$ is added to the graph. Then any path from $s$ to $t$ in the original graph, forms a cycle in the new graph.  It is known that all cycles belong to the kernel of the incidence matrix. The condition that the cycle belongs to the kernel of the new incidence matrix and contains the edge $(t,s)$  translates to the condition~\eqref{eq:Dxy} for the original graph. 
\end{remark}

Evidently, the conclusion~\eqref{eq:Dxy} does hold in the reverse direction: vectors that satisfy the linear constraint $D\xp=\yst$ may take fractional values and do not correspond to a valid incidence vector. For instance, any linear combination $x=ax^{(p_1)} + (1-a)x^{(p_2)}$ for $a\in [0,1]$ of the incidence vectors $x^{(p_1)}$ and $x^{(p_2)}$ of two distinct paths, $p_1$ and $p_2$, between $s$ and $t$, satisfies the constraint $Dx=\yst$. 
Yet, if the shortest path is unique, then a solution to~\eqref{eq:Dxy} with the least number of nonzero entries would necessarily correspond to this shortest path. Thus, although the exact equivalency does not hold for~\eqref{eq:Dxy}, the shortest path can be recovered from the ``sparsest'' solution to~\eqref{eq:Dxy}.

Now, a well known fact, that underlies techniques in modern compressive sensing \cite{tibshirani1996regression}, is that the $\ell_1$ norm can be used as a suitable surrogate for obtaining ``sparse'' solutions. Thus, we
propose as a relaxation to the shortest path problem the following:
\begin{equation}\label{eq:linear-prog}
{\rm arg}\min_{x\in \mathbb R^m}~\|Wx\|_1,\quad \text{s.t.}\quad Dx = \yst
\end{equation}
This is a linear programming problem and, for the reasons just noted, if the shortest path is unique, then the solution turns out to be integer-valued and corresponding to a valid incidence vector~\cite[Theorem 6.5 (Integrality Theorem), p.~186]{ahuja1988network}.

\subsection{Lasso formulation}
It is natural to consider the following $l_1$-regularized regression problem by replacing the constraint with a penalty term and changing variables by introducing $\beta=Wx$: 
\begin{equation}
\min_{\beta \in \mathbb R^m}~ \frac{1}{2}\|y-Q\beta\|_2^2 + \lambda\|\beta\|_1. \label{eq:lasso}
\end{equation}
Here $\lambda >0$ is the regularization parameter and $Q\triangleq DW^{-1}$. Problem~\eqref{eq:lasso} is known as the {\em lasso problem}. For $\lambda>0$, the solution of~\eqref{eq:lasso} is no longer equal to the shortest path.  However, in the limit as $\lambda \to 0$, the solution becomes exact.
Our motivation for exploring the formulation~\eqref{eq:lasso} is twofold:\\[.05in]
(i) As shown in Section~\ref{sec:LARS}, the LARS algorithm, designed to solve the lasso problem~\eqref{eq:lasso}, is equivalent to Dijkstra's algorithm, and\\[.05in]
(ii) as discussed in Section~\ref{sec:proximal}, it allows the flexibility to us proximal optimization methods to obtain a good approximation of the shortest path in large graph setting. 

%With a change of variable $\beta =Wx$, the  generalized lasso problem~\eqref{eq:lasso} becomes the standard lasso problem 
%\begin{equation}\label{eq:lasso-2}
%\min_{\beta \in \mathbb R^m}~ \frac{1}{2}\|y-\tilde{D}\beta\|_2^2 + \lambda\|\beta\|_1
%\end{equation}
%where we defined $\tilde{D}\triangleq DW^{-1}$ and dropped the superscript $(s,t)$ from $\yst$ for simplicity.  

%Alternatively, the KKT condition is expressed as
%\begin{equation*}
%\begin{cases}
%\lambda  = Q_j^T R\quad &\text{if}~\beta_j>0\\
%\lambda  = -Q_j^T R\quad &\text{if}~\beta_j<0\\
%\lambda  \geq  |Q_j^T R| \quad &\text{if}~\beta_j=0
%\end{cases}
%\quad \text{where} \ \ R = y-Q\beta 
%\end{equation*} 
%where $Q_j$ denotes the $j$-th column of the matrix $Q$.  
%A solution can be expressed as
%\begin{equation}
%\beta = ()
%\end{equation}

\subsection{Uniqueness of the Lasso solution}
The solution to the Lasso problem is unique when $\text{rank}(D)=m$, i.e. the rank of incidence matrix is equal to the number of edges. This condition holds if only if the graph is a tree (or a collection of disjoint  trees).  Evidently, the assumption that the graph is a tree is too restrictive, especially for the shortest path problem, because  the problem becomes trivial. 

The rank condition $\text{rank}(D)=m$ is sufficient but not necessary. Relaxations of this assumption have been introduced in the literature~\cite{tibshirani2013lasso}. We use the result~\cite[lemma 2]{tibshirani2013lasso} to prove the uniqueness of the lasso solution under the following assumption.

\noindent
{\bf Assumption A1:} The shortest path  between vertex $s$ or $t$ and any other vertex  is unique. 

The uniqueness result is expressed in the following lemma. The proof appears in Appendix~ \ref{apdx:proof-lemma-3.2}.\\[-.1in]

\begin{lemma}
	Under assumption A1, the lasso problem~\eqref{eq:lasso}  admits a unique solution for all $\lambda>0$.
\end{lemma}
%\AHT{Write the proof in appendix(added)}\\
%\begin{proof}
%	See details in Appendix \ref{apdx:proof-lemma-3.2}
%\end{proof}
\begin{remark}
	Assumption A1 is satisfied for a generic selection of weights, for example if a small noise is added to the weights. The assumption is necessary, as it is straightforward to come up with counter examples. 
%	Eventhough the uniqueness does not hold, it is still possible to simulate the LARS algorithm. See Section~\ref{sec:non-unqiue} to see how does LARS performs on a example with non-unique solution. 
\end{remark}

\section{Solution path and LARS algorithm}\label{sec:LARS}

%For a connected graph, This condition hold  under the assumption that the rank of $n\times m$ incidence  matrix $D$ is larger than $m$ . However, in the graph setting where $D$ is incidence matrix of a graph, the full-rank condition is not satisfied.   Let $\beta(\lambda)$ denote the solution to 
%% With a change of variable $z=W x$ and $D=DW^{-1}$ Let $x(\lambda)$ denote the minimizer of the 
%In this section, we derive the KKT condition of \eqref{eq:0} and provide the LARS algorithm

%The LARS algorithm  
%Since we want to use the LARS algorithm, let's first derive the KKT condition of the  lasso problem \eqref{eq:0}, this part is from the appendix of \cite{tibshirani2011solution}. 
%\begin{align}
%    \hat{\beta}(\lambda)=\min_{\beta} \frac{1}{2}\|Y - X\beta\|_{2}^{2} + \lambda |\beta |_{1}.\label{eq:12}  %%equation 12
%\end{align}
%Deriving the KKT condition of \eqref{eq:12} by using the fact that $\|\beta\|_{1}$ is subdifferentiable \cite{boyd2004convex} and the subgradient is
%%\begin{align}\label{eq:subgradient}
%%    \gamma_{i} \in
%%    \begin{cases}
%%      \text{sign}(\beta_{i}\} \ \ \  \text{if}\;\beta_{i} \neq 0 \\
%%      \left[ -1, 1 \right] \ \ \ \ \ \text{if}\; \beta_{i}=0
%%     \end{cases}  \ \ \text{for} \ \ i=1, ..., m.  
%%\end{align}
%We can get 
%\begin{align*}
%X^{T}(Y - X\beta) = \lambda\gamma 
%\end{align*}
\subsection{KKT conditions}
Let $\beta(\lambda)$ denote the solution to the lasso problem~\eqref{eq:lasso}. It must satisfy the  KKT condition,
%\subsection{KKT conditions}
%A vector $\beta \in \mathbb R^m$ minimizes the Lasso problem~\eqref{eq:lasso} if and only if the following 
%KKT conditions are satisfied
\begin{align}\label{eq:KKT}
Q^T(y - Q\beta(\lambda)) = \lambda\gamma,
\end{align}
where $\gamma$ belongs to  sub-differential of $\|\beta(\lambda)\|_1$ whose $j$-th component is given by
\begin{align}\label{eq:subgradient}
\gamma_{j} \in
\begin{cases}
\{\text{sign}(\beta_{j}(\lambda))\}~~  &\text{if}\;\beta_{j}(\lambda) \neq 0 \\
\left[ -1, 1 \right] ~~ &\text{if}\; \beta_{j}(\lambda)=0.
\end{cases}  
\end{align}
%and we dropped the superscript $(s,t)$ of $\yst$ for simplicity. 

The KKT condition motivates to divide the indices $\{1,\ 2,\dots, m\}$ into two sets: active set $\calA$, where $\beta(\lambda)$ is nonzero, and non-active set $\calA^c$, where $\beta(\lambda)$ is zero. 
Let $\beta_\calA(\lambda)$ denote the vector $\beta(\lambda)$ where non-active components are removed, and $Q_\calA$ be the matrix $Q$ where the columns corresponding to non-active set are removed. Then, the KKT condition~\eqref{eq:KKT} is expressed as 
\begin{subequations}\label{eq:KKT-new}
	\begin{align}\label{eq:KKT-active}
	Q_j^T(y-Q_\calA\beta_\calA(\lambda)) &= s_j \lambda,\quad \forall j \in \calA\\\label{eq:KKT-non-active}
	Q_j^T(y-Q_\calA\beta_\calA(\lambda)) &\in [-\lambda,\lambda],\quad \forall j \in \calA^c,
	\end{align}
\end{subequations}
where $Q_j$ denotes the $j$-th column of $Q$ and the sign vector
\begin{equation}
s \triangleq \text{sign}\left(Q_\calA^{T}(y - Q_\calA\beta_\calA(\lambda )) \right) = \text{sign}\left(\beta_\calA \right).
\end{equation}

%The active set $\calA$, is the the subset of indices where  $\beta_j(\lambda)\neq 0$ and the non-active set $\calA^c$ is is the subset of indices where $\beta_j(\lambda)= 0$. 
%On the active set, the KKT condition~\eqref{eq:KKT} is satisfied with $\gamma_j = \pm 1$, i.e. 
%\begin{align}\label{eq:active-set}
%Q_{\calA}^{T}(y - Q_\calA\beta_\calA(\lambda )) = \lambda s %\quad \forall j \in \calA
%\end{align}
%where $Q_\calA$ stands for the columns of the matrix $Q$ that belong to the active set, $\beta_\calA$ is the components of $\beta$ in the active set,  and
%\begin{equation}
%s \triangleq \text{sign}\left(Q_\calA^{T}(y - Q_\calA\beta_\calA(\lambda )) \right) = \text{sign}\left(\beta_\calA \right)
%\end{equation}
%
%On the non-active set, the KKT condition~\eqref{eq:KKT} is satisfied with $\gamma_j \in [-1,1]$, i.e.
%\begin{align}\label{eq:non-active-set}
%Q_{j}^{T}(y - Q_\calA\beta_\calA(\lambda )) 
%\in [-\lambda,\lambda], \quad \forall j \in \calA^c
%\end{align}
%The two conditions~\eqref{eq:active-set}-\eqref{eq:non-active-set} are summarized together as
%\begin{subequations}\label{eq:KKT-new}
%\begin{align}\label{eq:KKT-active}
%Q_j^T(y-Q_\calA\beta_\calA(\lambda)) &= s_j \lambda,\quad \forall j \in \calA\\\label{eq:KKT-non-active}
%Q_j^T(y-Q_\calA\beta_\calA(\lambda)) &\in [-\lambda,\lambda],\quad \forall j \in \calA^c
%\end{align}
%\end{subequations}

\subsection{LARS algorithm}
The least angle regression (LARS) algorithm,
in its lasso state\footnote{Whereas the original LARS algorithm does not provide the lasso solution, a modification in~\cite{efron2004least} does indeed solve the lasso problem.},
finds the solution $\beta(\lambda)$ that satisfy the KKT condition~\eqref{eq:KKT-new} for all $\lambda>0$~\cite{efron2004least}.
The vector $\beta(\lambda)$ is continuous and piecewise linear, as a function of $\lambda$, with break points  $\lambda_1>\lambda_2>\ldots>\lambda_l>0$. The active set and the sign vector remain the same during each interval and change at each break points.  Let $\calA_k$ and $s_k$ denote the active set and sign vector during the  interval $(\lambda_{k+1},\lambda_{k})$. The LARS algorithm starts with $\lambda_0=\infty$, $\calA_0=\emptyset$, and $s_0=\emptyset$. Then, at iteration $k$, given $\lambda_k$, $\calA_k$, and $s_k$, the algorithm finds the next breaking point $\lambda_{k+1}$, the next active set $\calA_{k+1}$, and the next sign vector $s_{k+1}$. During each interval $(\lambda_{k+1},\lambda_{k})$, the  vector $\beta(\lambda)$ is the minimum $l_2$-norm solution of the condition~\eqref{eq:KKT-active} given by:
\begin{align}\label{eq:betaA}
\beta_{\calA_k}(\lambda) &= 
(Q_{\calA_k}^{T}Q_{\mathcal{A}_k})^{+}(Q_{\mathcal{A}_k}^{T}y-\lambda s_k)\\\nonumber
&= a^{(k)} - b^{(k)} \lambda , %+ r
\end{align} 
where $^+$ denotes the Moore-Penrose pseudo-inverse and 
\begin{subequations}\label{eq:a-b}
	\begin{align}
	a^{(k)} &\triangleq (Q_{\calA_k}^{T}Q_{\mathcal{A}_k})^{+}Q_{\mathcal{A}_k}^{T}y,\\ b^{(k)}&\triangleq (Q_{\calA_k}^{T}Q_{\mathcal{A}_k})^{+}s_k.
	\end{align}
\end{subequations}
The next breaking point $\lambda_{k+1}$ is the largest value $\lambda<\lambda_k$ so that \eqref{eq:betaA} does not satisfy the KKT conditions~\eqref{eq:KKT-new} anymore. The KKT conditions are violated in two cases: 
\begin{itemize}
	\item Joining: This case happens when condition~\eqref{eq:KKT-non-active} is violated for some $j\in \calA^c_k$, i.e.   $|Q_j^T(y-Q_{\calA_k}\beta_{\calA_k}(\lambda))| = \lambda$. 
%	There is an index $j\in \calA^c$ so that $|Q_j^T(y-Q_\calA\beta_{\calA_k}(\lambda))| = \lambda$ and condition~\eqref{eq:KKT-non-active} does not hold anymore. 
	For each index $j\in \calA^c_k$, this happens at $\lambda= t^{\text{join}}_j$ given by
	\begin{equation}\label{eq:join}
	t^{\text{join}}_{j,k} = \frac{Q_j^T(Q_{\calA_k}a^{(k)}-y)}{Q_j^TQ_{\calA_k}b^{(k)} \pm 1}.
	\end{equation}
	\item Crossing: This case happens when condition~\eqref{eq:KKT-active} is violated for some $j\in \calA_k$. By definition of $\beta_{\calA_k}(\lambda)$ according to~\eqref{eq:betaA}, the condition~\eqref{eq:KKT-active}  is violated only when $s\neq s_k$. This happens when one of the component of $\beta_{\calA_k}(\lambda)$ crosses zero, i.e. $a^{(k)}_j - \lambda b^{(k)}_j=0$ for some $\lambda<\lambda_k$. For each index $j\in \calA_k$, the crossing happens at $\lambda=t^{\text{cross}}_j$ given by 
	\begin{equation}\label{eq:cross}
	t_{j,k}^{\text{cross}} = \begin{cases}
	\frac{a^{(k)}_j}{b^{(k)}_j}\quad &\text{if}~ 0 < \frac{a^{(k)}_j}{b^{(k)}_j} < \lambda_k\\
	0 \quad &\text{otherwise}.
	\end{cases}
	\end{equation} 	
\end{itemize}  
The algorithm selects the next break point $\lambda_{k+1}$ to be the maximum of joining times and crossing times:
\begin{equation}\label{eq:next-break}
\lambda_{k+1}=\max (\max_{j \in \calA_k^c} 	t^{\text{join}}_{j,k}, \max_{j \in \calA_k} 	t_{j,k}^{\text{cross}} )
\end{equation}
If the joining happens, the joining index is added to the active set and the sign vector is updated accordingly. If a crossing happens, the crossing index is removed from the active set. The overall algorithm is summarized in~\ref{alg:LARS}.

\begin{algorithm}[H]
	\caption{\bf LARS path algorithm for the lasso problem~\eqref{eq:lasso}}
	\begin{algorithmic}[1]
		\renewcommand{\algorithmicrequire}{\textbf{Input:}}
		\renewcommand{\algorithmicensure}{\textbf{Output:}}
		\REQUIRE matrix $Q=DW^{-1}$ and vector $y=\yst$
		\ENSURE  incidence vector $x=W^{-1}\beta$ and  path length $\|\beta\|_1$
		\STATE  $k=0$, $\lambda_{0}=\infty$, $\mathcal{A}=\emptyset$, $s=0$, $a^{(0)} =0$ and $b^{(0)}=0$. 
		\WHILE{$\lambda_{k}>0$} 	
%		 by least squares
%		\begin{align*}
%		\hat{\beta}_{\lambda_{k},\mathcal{A}} = ((D_{\mathcal{A}})^{T}D_{\mathcal{A}})^{+}((D_{\mathcal{A}})^{T}(-Y_{st})-\lambda_{k}s).
%		\end{align*}
		\STATE Compute the joining time~\eqref{eq:join} for $j\in \calA^c_k$.
		\STATE Compute the crossing time~\eqref{eq:cross} for $j\in \calA_k$.
		\STATE Set $\lambda_{k+1}$ according to~\eqref{eq:next-break} 
		\begin{itemize}
			\item If join happens, add the joining index to $\mathcal{A}_k$ and its sign to $s_k$.
			\item If cross happens, remove the crossing index from $\mathcal{A}_k$ and its sign from $s_k$.
		\end{itemize}
		\STATE $k = k+1$.
		\STATE Compute $a^{(k)}$ and $b^{(k)}$ according to~\eqref{eq:a-b}
		\STATE Set $\beta_{\calA_k}  = a^{(k)} - \lambda_kb^{(k)}$ and $\beta_{\calA_k^c}=0$.
		\ENDWHILE
		\RETURN $x=W^{-1}\beta$ and $\|\beta\|_1$
	\end{algorithmic}
\label{alg:LARS}
\end{algorithm}

\subsection{Numerical example}
Consider the Nicholson's graph \cite[p.~6]{pohl1969bidirectional}, as depicted in Figure~\ref{fig:a}, and the shortest path problem between vertex $1$ and vertex $9$. The iterations of the LARS algorithm are depicted in Figure~\ref{fig:example}. At each iteration, the edges that belong to the active set $\mathcal{A}$ are highlighted in red. It is observed that at each iterations, edges are added to the active set and are never removed. The algorithm terminates after four iterations when $\lambda_5=0$ and a path between vertex $1$ and $9$ is formed. The lasso solution path $\beta(\lambda)$ is depicted in Figure~\ref{fig:lassopath}.

Example~\ref{fig:example} illustrates the similarity between LARS algorithm and  the Dijkstra's algorithm. Namely, the LARS algorithm constructs two shortest-path trees, with roots at vertex $1$ and vertex $9$ respectively. This is similar to the bi-directional  Dijkstra's algorithm, as discussed in Section~\ref{sec:prelim}. In the next section, we show that the similarity between the LARS algorithm and the Dijkstra's algorithm holds in general. 

%This is similar to the bi-directional Dijkstra's algorithm. It is known that the Dijkstra's algorithm also constructs the shorest   

%  $\mathcal{T}_{s}$ and $\mathcal{T}_{t}$ treat $v_{s}$ and $v_{t}$ to be the roots. That is equal to say that the algorithm never remove elements from $\mathcal{A}$, just adding. More than that, the spanning tree is actually the shortest path tree. We will prove the above points in the next section.

\section{Relationship between  LASSO and Dijkstra}\label{sec:relation}
We establish the connection between the LARS algorithm~\ref{alg:LARS} and the Dijkstra's algorithm by showing that the LARS algorithm iteratively builds two shortest path trees with roots at $s$ and $t$, and that the algorithm terminates when the two trees connect.

We prove this by induction. The induction hypothesis is as follows: At iteration $k$ of the algorithm, the edges in the active set $\calA_k = \calA_k^{(s)} \cup \calA_k^{(t)}$ form two disjoint subsets $\calA_k^{(s)}$ and $\calA_k^{(s)}$.   Each subset form a tree on the vertices,
%induction assumption is as follows: The active set $\calA_k= \calA^{(s)}_k \cup \calA^{(t)}_k$ is comprised of two set of disjoint edges $\calA^{(s)}_k$ and $\calA^{(t)}_k$. The two set of edges build two 
denoted by $T^{(s)}_k \subset \calV$ with root at $s$ and $T^{(t)}_k \subset \calV$ with  root at $t$, respectively. The two trees are the shortest-path trees from the root vertex. Moreover, crossing does not occur at this iteration, i.e., no edges are removed from the active set. 

The induction hypothesis is true at $k=0$,  because the active set is empty,  the two trees consist of single root vertex, i.e.  $T^{(s)}_0=\{s\}$ and $T^{(t)}_0=\{t\}$, and crossing does not occur because the active set is empty. 

Assuming the induction hypothesis at iteration $k$, we show  the hypothesis also holds at iteration $k+1$ by proving:
\begin{itemize}
	\item[(a)] Let $v_\text{min}^{(s)}$ and $v_\text{min}^{(t)}$  denote the vertex that has the minimum distance to the root $s$ and $t$ respectively, among all vertices outside the two trees. Then, either the edge that  connects vertex $v_\text{min}^{(s)}$ to tree $T^{(s)}_k$ or the edge that connects $v_\text{min}^{(t)}$ to $T^{(t)}_k$ is added to the active set; 
	
	\item[(b)] Crossing does not occur. 
	%	Edges are only added.
\end{itemize}
%Note that (a) proves that  the two trees remain shortest path  trees at iteration $k+1$.  

%%{\color{red}Also, they can be denoted by $\mathcal{A}_{k}^{(s)} \subset \mathcal{E}$ and $\mathcal{A}_{k}^{(t)} \subset \mathcal{E}$.}
%The hypothesis is true at i
%We need to show that at the next iteration, the new active-set $\calA_{k+1}$ also gives two shortest path trees $T^{(s)}_{k+1}$  and $T^{(t)}_{k+1}$ rooted at $s$ and $t$ respectively. In particular, we need to show:
%\begin{itemize}
%	\item[(a)] Let $v_\text{min}^{(s)}$ and $v_\text{min}^{(t)}$  denote the vertex that has the minimum distance to the root $s$ and $t$ respectively, among all vertices outside the two trees. Then, either the edge that  connects vertex $v_\text{min}^{(s)}$ to tree $T^{(s)}_k$ or the edge that connects $v_\text{min}^{(t)}$ to $T^{(t)}_k$ is added to the active set;
%	   
%	\item[(b)] Crossing does not occur, hence edges are not removed from the active set. Edges are only added.
%\end{itemize}
Moreover, we also need to show the termination condition
\begin{itemize}
	\item[(c)] The algorithm terminates when the two trees connect. 
\end{itemize}
%Note that the induction hypothesis is true at $k=0$ by defining $T^{(s)}_0=\{s\}$ and $T^{(t)}_0=\{t\}$.  Next, we prove (a), (b) and (c).  

The proof is based on simplified expressions for joining time and the crossing time that are obtained using Lemma~\ref{lem:incidence}. The derivations appear in Appendix~\ref{apdx:t-cross-derivation} and~\ref{apdx:t-join-derivation}. %The result for the joining time is stated first. 
%We divide the edges that do not belong to the active set into five groups:
%\begin{itemize}
%	\item $\calT^2$: edges connecting two vertices in the same tree
%	\item $\Omega^2$: edges connecting two vertices outside the two trees
%\item $\calT^{(s)} \times \Omega$: edges connecting a vertex from tree $T^{(s)}_k$ to an outside vertex   
%\item $\calT^{(t)} \times \Omega$: edges connecting a vertex from tree $T^{(t)}_k$ to an outside vertex   
%\item $\calT^{(s)} \times \calT^{(t)}$: edges connecting  a vertex from $T^{(s)}_k$  to a vertex from $T^{(t)}_k$
%\end{itemize}

\subsubsection{Joining time} 
For the edge $e_j=(v_1,v_2)$, where $e_j\in \calA_k^c$, the joining time is  
\begin{align} \label{eq:tjoin}
t^{\text{join}}_{j,k} = 
\begin{cases}
0 ~&\text{if}~(v_1,v_2) \in \Omega^2   \\
0 ~&\text{if}~(v_1,v_2) \in {T^{(s)}_k}^2 \cup {T^{(t)}_k}^2    \\
\frac{1}{|T^{(s)}_k|l^{(s)}_{v_2}  - \sum_{v\in T^{(s)}_k} l^{(s)}_v} ~ &\text{if}~(v_1,v_2) \in {T^{(s)}_k} \times \Omega \\
\frac{1}{|T^{(t)}_k|l^{(s)}_{v_2}  - \sum_{v\in T^{(t)}_k} l^{(t)}_v } ~ &\text{if}~(v_1,v_2) \in {T^{(t)}_k} \times \Omega\\
\frac{|T^{(s)}_k|+|T^{(t)}_k|}{\gamma  } ~ &\text{if}~(v_1,v_2) \in {T^{(s)}_k} \times{T^{(t)}_k}  
\end{cases}
\end{align}
where $\Omega$ is the set of vertices not in the trees, $l^{(s)}_v$ and $l^{(t)}_v$ denote the distance of vertex $v$ to the root $s$ and $t$ respectively, and \[\gamma =|\Ts_k| |\Tt_k| l^{(s)}_t - |\Tt_k| \sum_{v\in \Ts_k}l_v^{(s)} - |\Ts_k|\sum_{v\in \Tt_k}l_v^{(t)}.\]   
% For a vertex $v$ in the tree $T^{(s)}_k$ or $T^{(s)}_k$, let $l^{(s)}_v$ denote the distance of the vertex $v$ to the root $s$. Similarly, define $l^{(t)}_v$ for vertices $v \in T^{(t)}_k$. 
%\begin{subequations}
%	\begin{align}
%	\hat{l}^{(s)} &\triangleq \frac{1}{|T^{(s)}_k|}\sum_{v\in T^{(s)}_k} l^{(s)}_v,~~
%	\hat{l}^{(t)} \triangleq \frac{1}{|T^{(t)}_k|}\sum_{v\in T^{(t)}_k} l^{(t)}_v
%	\end{align}
%\end{subequations} are the average distance of the vertices in a tree. 
%The details of the derivation appears in Appendix~\ref{apdx:t-cross-derivation} \AHT{add the reference to appendix(Done)}
\subsubsection{Crossing time} 
For an edge $e_j=(v_1,v_2)$ where $e_j\in \calA_k$, the expression $a_{j}^{(k)}/b_j^{(k)}$ that appears in the definition of crossing time~\eqref{eq:cross} is
\begin{align}\label{eq:crossingtime}
\frac{a_j^{(k)}}{b_j^{(k)}}= 
\begin{cases}
\frac{1}{\frac{|\Ts_k|}{|R^{(s)}_j|} \sum\limits_{v \in R^{(s)}_j} l^{(s)}_v - \sum\limits _{v\in T^{(s)}_k} l^{(s)}_v} &~~\text{if}~(v_1,v_2)\in {T^{(s)}_k}^2 \\
\frac{1}{\frac{|\Tt_k|}{|R^{(t)}_j|} \sum\limits _{v \in R^{(t)}_j} l^{(t)}_v - \sum\limits _{v\in T^{(t)}_k} l^{(t)}_v} &~~\text{if}~(v_1,v_2)\in {T^{(t)}_k}^2 
\end{cases}
\end{align}
where $R^{(s)}_j$ and $R^{(t)}_j$  are the subsets of vertices in the tree $T^{(s)}_k$ and $T^{(t)}_k$  respectively, whose path to the root contains the edge $e_j$. 
%The details of the derivation appears in Appendix~\ref{apdx:t-join-derivation}\\[0.1in]

Proof of (a): Assume there is no edge that connects the two trees. i.e. the last case in expression~\eqref{eq:tjoin} does not happen. We study this case in part (c). Then, the maximum of $t^{\text{join}}_{j,k}$ is given by 
%Using the simplified expression for the joining time~\eqref{eq:tjoin}
\begin{align*}
  \max (\frac{1}{|\Ts_k|l_{v^{(s)}_{\min}}  -\sum\limits _{v\in T^{(s)}_k} l^{(s)}_v }, \frac{1}{|\Tt_k|l_{v^{(t)}_{\min}}  - \sum\limits _{v\in T^{(s)}_k} l^{(s)}_v})
\end{align*}
where  the first expression is achieved by the edge that connects $v_\text{min}^{(s)}$  to three $T^{(s)}_k$ and the second expression is achieved by the edge that connects $v_\text{min}^{(t)}$  to tree $T^{(t)}_k$. Hence, one of these two edges is joined to the active set, if crossing does not occur. In part-(b), we show crossing does not occur.  

\medskip
\begin{remark}[No cycles] Cycles may created in the following two scenarios: (i) an edge that connects two vertices of a tree is joined; (ii) Two edges that connect the tree to a single vertex, say $v$, are joined simultaneously. The scenario (i) can not happen because $t_j^{\text{join}}=0$ for such edges (second case in~\eqref{eq:tjoin}). The scenario (ii) can not happen, because in order for two edges to join simultaneously, we must have two  distinct shortest path from $v$ to the root, which is not possible according to Assumption A1. 
\end{remark}
%\begin{align*}
%l_{v^{(s)}_{\min}} = \min (l_{v^{(s)}_{1}} + w_{j}).\\
%\max_{e_{j} \in \mathcal{A}_{k}^{(t)}}t^{ join}_{j,k} = 
%\frac{1/n_t}{l_{v^{(t)}_{\min}}  - \hat{l}^{(t)} }~~\text{s.t.}~l_{v^{(t)}_{\min}} = \min (l_{v^{(t)}_{1}} + w_{j}). 
%\end{align*}
%Before $T_{k}^{(s)}$ and $T_{k}^{(t)}$ are connect, $t_{j}^{join} \neq 0,~~\forall~
%e_{j} \in \mathcal{A}_{k}$. In the $k$-th iteration, $n_{s}$, $n_{t}$, $\hat{l}^{(s)}$ and $\hat{l}^{(t)}$ are fixed, if the joining happens in $k$-th iteration:
%\begin{align*}
%\max_{e_{j} \in \mathcal{A}_{k}^{(s)}}t^{ join}_{j,k} = 
%\frac{1/n_s}{l_{v^{(s)}_{\min}}  - \hat{l}^{(s)} }~~\text{s.t.}~l_{v^{(s)}_{\min}} = \min (l_{v^{(s)}_{1}} + w_{j}).\\
%\max_{e_{j} \in \mathcal{A}_{k}^{(t)}}t^{ join}_{j,k} = 
%\frac{1/n_t}{l_{v^{(t)}_{\min}}  - \hat{l}^{(t)} }~~\text{s.t.}~l_{v^{(t)}_{\min}} = \min (l_{v^{(t)}_{1}} + w_{j}). 
%\end{align*}
%\begin{align*}
%%    t^{join}_{k} = 
%    \max_{e_j \in \mathcal{A}_{k}} t^{ join}_{j,k} = 
%    \max(\max_{e_j \in \mathcal{A}_{k}^{(s)}}t^{ join}_{j,k}, \max_{e_j \in \mathcal{A}_{k}^{(t)}}t^{ join}_{j,k})
%\end{align*}
%and 
%
%Thus,
%\begin{align}
%    t^{join}_{k} = \max_{j \in \mathcal{A}_{k}} (\frac{1/n_s}{l_{v^{(s)}_{\min}}  - \hat{l}^{(s)}}, \frac{1/n_t}{l_{v^{(t)}_{\min}}  - \hat{l}^{(t)}})
%\end{align}
%Therefore,  the edge that connects vertex $v_\text{min}^{(s)}$ to tree $T^{(s)}_k$ or the edge that connects $v_\text{min}^{(t)}$ to $T^{(t)}_k$ is  added to the active set $\mathcal{A}_{k+1}$. 

Proof of (b): To prove crossing does not occur, we show that $a^{(k)}_j/b^{(k)}_j\geq \lambda_{k}$ for all $e_j$ in the active set, so that crossing time is zero according to the definition~\eqref{eq:cross}. In order to do so, first we obtain an expression for $\lambda_k$ and then compare it to crossing times. $\lambda_k$ is determined by the maximum of joining time and crossing time at iteration $k+1$ according to~\eqref{eq:next-break}. By induction assumption, crossing did not occur in the iteration $k-1$. Hence, $\lambda_k$ is determined by the maximum joining time, which by part-(a) takes two possible values, corresponding to the edge that connects to tree $T^{(s)}_{k-1}$ or the edge that connects to  tree $T^{(t)}_{k-1}$. Without loss of generality, assume the joining happens to the tree $T^{(s)}_{k-1}$. Then, 
\begin{align*}
\lambda_{k-1}&=  \frac{1}{|\Ts_{k-1}|l_{v^{(s)}_{\min}}  -\sum_{v\in T^{(s)}_{k-1}} l^{(s)}_v} 
%\\
%&=  \frac{1}{|T^{(s)}_{k-1}|l_{v^{(s)}_{\min}}  - \sum_{v\in T^{(s)}_{k-1}} l^{(s)}_v}
\end{align*}   
Next, we show $a^{(k)}_j/b^{(k)}_j\geq \lambda_{k-1}$ for all $e_j$ that belong to the tree $T^{(s)}_k$. For such an edge we have from~\eqref{eq:crossingtime} that 
\begin{align*}
\frac{a_j^{(k)}}{b_j^{(k)}}
%&= 
%\frac{\frac{1}{|T^{(s)}_{k}|}}{\frac{1}{|R^{(s)}_j|} \sum_{v \in R^{(s)}_j} l^{(s)}_v -  \frac{1}{|T^{(s)}_{k}|}\sum_{v\in T^{(s)}_{k}} l^{(s)}_v} \\
&=\frac{1}{\frac{1+|T^{(s)}_{k-1}|}{|R^{(s)}_j|} \sum_{v \in R^{(s)}_j} l^{(s)}_v - l_{v^{(s)}_{\min}} - \sum_{v\in T^{(s)}_{k-1}} l^{(s)}_v}\\
&\geq \frac{1}{|T^{(s)}_{k-1}| l_{v^{(s)}_{\min}} - \sum_{v\in T^{(s)}_{k-1}} l^{(s)}_v}
\end{align*}
where we used $|T^{(s)}_k|=|T^{(s)}_{k-1}|+1$, $\sum_{v\in T^{(s)}_{k}} l^{(s)}_v=  l_{v^{(s)}_{\min}}+ \sum_{v\in T^{(s)}_{k-1}} l^{(s)}_v$, and $l_{v^{(s)}_{\min}}  \geq l_{v}$ for all $v\in T^{(s)}_k$. The last statement is true because $v^{(s)}_{\min}$ is the latest vertex that is added to the tree and other vertices that have been already added have a shorter distance to the root. 

The proof that $a^{(k)}_j/b^{(k)}_j\geq \lambda_{k-1}$ for all $e_j$ that belong to the other tree $T^{(t)}_k$ is conceptually similar. One needs to compare $a^{(k)}_j/b^{(k)}_j$ with the joining time of the last edge that has been added to the tree $T^{(t)}_k$ at a certain past iteration, say $k'<k$, and use the fact that $\lambda_k<\lambda_{k'}$. The details are omitted on account of space. 

\begin{figure}[htb!]
  	\centering
  	\scalebox{1}{
  		\begin{tikzpicture}
  		[process/.style={circle,draw=blue!50,fill=blue!20,thick,
  			inner sep=0pt,minimum size=7mm},
  		target/.style={rectangle,draw=black!50,fill=black!20,thick,
  			inner sep=0pt,minimum size=8mm}]
  		\node (n1) at (0,0)   [process] {$v_{s}$};
  		\node (n2) at (1.5,0) [process] {$v_{1}$};
  		\node (n3) at (3.75,0) [process] {$v_{2}$};
  		\node (n4) at (6.25,0) [process] {$v_{3}$};
  		\node (n5) at (8.5,0) [process] {$v_{4}$};
  		\node (n6) at (10,0)  [process] {$v_{t}$};
  		
  		\draw[dashed][ultra thick] (n1) to   (n2);
  		\draw[ultra thick] (n2) to  node[midway,above](){$e_{p}$} (n3);
  		\draw[ultra thick] (n3) to  node[midway,above](){$e_{i}$} (n4);
  		\draw[ultra thick] (n4) to  node[midway,above](){$e_{s}$} (n5);
  		\draw[dashed][ultra thick] (n5) to  (n6);
  		\end{tikzpicture}}
  	\caption{Path $\mathcal{P}_{s,t}$}
  	\label{fig:part-c}
  \end{figure}
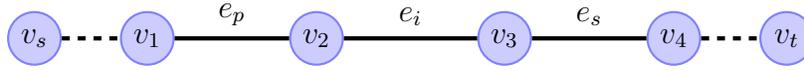

Proof of (c): Assume  the two trees $\Ts_k$ and $\Tt_k$ become connected at iteration $k$. This happens when the last expression in~\eqref{eq:tjoin} achieves the maximum joining time, hence
\begin{align*}
\lambda_k = \frac{|\Ts_k| + |\Tt_k|}{\gamma}.
\end{align*}
The situation is depicted in Figure~\ref{fig:part-c}. The objective is to show that the algorithm terminates after this, i.e. $\lambda_{k+1}=0$. We show this by proving the joining time and crossing time are both zero. The derivation of joining time in Appendix~\ref{apdx:t-join-derivation} reveals that  
\begin{align*}
t_{j,k+1}^{join} = 0, &~~\forall e_j \in \mathcal{A}^{c}_{k+1}.
\end{align*}
For the crossing time, the derivation of Appendix~\ref{apdx:t-cross-derivation} yields that for all $e_j \in \calA_{k+1}$,
\begin{align}\label{eq:tcrossonetree}
\frac{a_j^{(k+1)}}{b_j^{(k+1)}} =
%& = \frac{D_{\calA^{+}_{k+1}}y^{(s,t)}}{D^{+}_{\calA_{k+1}}(D^{+}_{\calA_{k+1}})^{T}Ws} \nonumber \\
&\begin{cases}
0,~ \text{if}~e_j \notin p_{s,t}\\
%\frac{1}{ \gamma  + \frac{1}{|T^{(s)}_k|+|T^{(t)}_k|} (|\Ts_k|\sum_{v\in R_j^{(s)}}l_v^{(s)} - |\Tt_k| \sum_{v\in R_j^{(t)}} l_v^{(t)})} 
\frac{1}{\sum_{v\in R_j} l_v^{(s)} - \frac{|R_j|}{|\Tt_k|+|\Ts_k|} \sum_{v \in \Ts_k \cup \Tt_k} l_v^{(s)} },~\text{else},
\end{cases} 
\end{align} 
where $p_{s,t} \subset \calA_{k+1}$ is the path from $s$ to $t$. Therefore, it remains to show that the crossing time for the edges in $\calA_{k+1}\cap p_{s,t}$ are zero. We show this by proving $a_j^{(k+1)}/b_j^{(k+1)}\geq \lambda _k$ for all edges $e_j\in \calA_{k+1}\cap p_{s,t}$. First consider the edges that belong to the tree $\calA^{(s)}_{k}$. Then for these edges we have 
\begin{align*}
\frac{|\Ts_k|+|\Tt_k|}{\frac{a_j^{(k+1)}}{b_j^{(k+1)}}} =\gamma  + |\Ts_k| \sum_{v\in R_j^{(s)}}l_v^{(s)}- |R_j^{(s)}|\sum_{v \in \Ts_k} l_v^{(s)} 
 - |\Tt_k| \sum_{v\in R_j^{(s)}} l_v^{(t)} + |R_j^{(s)}|\sum_{v \in \Tt_k} l_v^{(t)} 
\end{align*}
where $R_j^{(s)}$ are the vertices in the tree $\Ts_k$ such that their path to the root $s$ contains $e_j$.  Because $l_v^{(s)} \leq l_{v_2}^{(s)}$ and $l_v^{(t)} \geq l_{v_2}^{(t)}$  for all $v\in R_j^{(s)}$, we have the inequality 
\begin{align*}
\frac{|\Ts_k|+|\Tt_k|}{\frac{a_j^{(k+1)}}{b_j^{(k+1)}}} -\gamma \leq |R_j|\bigg( |\Ts_k| l_{v_2}^{(s)}-\sum_{v \in \Ts_k} l_v^{(s)}
 - |\Tt_k|  l_{v_2}^{(t)} + \sum_{v \in \Tt_k} l_v^{(t)} \bigg).
\end{align*}
We claim that the expression in parentheses is negative
\begin{equation}\label{eq:claim}
|\Ts_k| l_{v_2}^{(s)}-\sum_{v \in \Ts_k} l_v^{(s)}  - |\Tt_k|  l_{v_2}^{(t)} + \sum_{v \in \Tt_k} l_v^{(t)} \leq 0.
\end{equation} 
If the claim is true, then
\begin{align*}
\frac{a_j^{(k+1)}}{b_j^{(k+1)}} \geq  \frac{|\Ts_k| + |\Tt_k|}{\gamma} = \lambda_k,
\end{align*}
proving that the crossing time is zero for edges $e_j \in \calA^{(s)}_k$.

Now we prove the claim. The edges are added in the order $e_s,e_p,e_i$ or $e_p,e_s,e_i$. In the first case, the joining time for the edges $e_p$ and $e_i$ are:   
\begin{align*}
t_{p,k-1}^{\text{join}} &= \frac{1}{|\Ts_k|l^{(s)}_{v_{2}} 
	- \sum\limits_{v\in T^{(s)}_k} l^{(s)}_v }\\
t_{i,k-1}^{\text{join}} &= \frac{1}{ |\Tt_k|l^{(t)}_{v_{2}} 
	-\sum\limits_{v\in T^{(t)}_k} l^{(t)}_v}.
\end{align*}  
The assumption that $e_p$ is added before $e_i$ implies $t_{p,k-1}^{\text{join}} >t_{i,k-1}^{\text{join}} $ concluding the claim~\eqref{eq:claim}. In the second case,  the joining time for the edges $e_p$ and $e_s$ are:   
\begin{align*}
t_{p,k-2}^{\text{join}} &= \frac{1}{|\Ts_k|l^{(s)}_{v_{2}} 
	- \sum\limits_{v\in T^{(s)}_k} l^{(s)}_v }\\
t_{s,k-2}^{\text{join}} &= \frac{1}{ |\Tt_k|l^{(t)}_{v_{2}} 
	-\sum\limits_{v\in T^{(t)}_k} l^{(t)}_v}.
\end{align*}
The order $e_p$ is added before $e_s$ concludes the claim~\eqref{eq:claim} because $t_{p,k-1}^{\text{join}} >t_{s,k-1}^{\text{join}}$. 

The proof that the crossing times for the edges that belong to the tree $\calA_k^{(t)}$ is by symmetry and interchanging $s$ and $t$.

\section{Proximal algorithm for large scale graph}\label{sec:proximal}

The alternating direction method of multipliers (ADMM) is a numerical algorithm that is used to solve a wide range of large-scale convex optimization problems~\cite{boyd2011distributed}.  Application of the ADMM to the lasso problem, as presented in~\cite[Section 6.4]{boyd2011distributed}, is based on the reformulation of the lasso  problem~\eqref{eq:lasso} as follows: 
\begin{equation}\label{eq:lasso-new}
\min_{\beta, \alpha \in \mathbb R^m}~ \frac{1}{2}\|y-Q\beta\|_2^2 + \lambda\|\alpha\|_1 +\frac{\rho}{2}\|\beta-\alpha\|^{2}_{2} ,\quad \text{s.t}\quad \alpha = \beta,
\end{equation}
where $\alpha \in  \mathbb{R}^m$ is an additional optimization variable, and $\rho>0$ is a positive constant. The Lagrangian corresponding to the constrained optimization problem~\eqref{eq:lasso-new} is 
\begin{align*}
L_{\rho}(\beta,\alpha,u) = \frac{1}{2}\|y-Q\beta\|_{2}^{2}+\lambda\|\alpha\|_{1} + u^{T}(\beta-\alpha)+\frac{\rho}{2}\|\beta-\alpha\|^{2}_{2}
\end{align*}
where $u\in \mathbb R^m$ is the Lagrange multiplier. Let $v := u/\rho$. The ADMM algorithm computes the optimal variables $\alpha,\beta,v$ iteratively according to
%and the ADMM steps, denote $w = u/\rho$ :
\begin{align}
\beta^{k+1} &= (Q^{T}Q+\rho I)^{-1}(Q^{T}y+\rho(\alpha^{k}-v^{k})) \label{eq:beta-update}\\
\alpha^{k+1} &= S_{\lambda/\rho}(\beta^{k+1}+\frac{1}{\rho}v^{k}) \nonumber\\
v^{k+1} &= v^{k}+\rho(\beta^{k+1} - \alpha^{k+1}) \nonumber
\end{align}
where $k$ is the iteration number,  and  $S_{\lambda/\rho}$ is the  soft-thresolding operator.  

The computational complexity of the ADMM iterations is dominated by the matrix inversion $(Q^{T}Q+\rho I)^{-1}$, which is of order $\mathcal{O}(p^3)$ (e.g. with \textit{Cholesky} decomposition),  where $p$ is  the number of edges. The complexity can be reduced using the matrix identity
\begin{equation}\label{eq:matrix-identity}
(Q^TQ + \rho I)^{-1} = \frac{1}{\rho}(I - Q^T(QQ^T + \rho I)^{-1}Q). 
\end{equation} 
which instead involves the matrix inversion $(QQ^T + \rho I)^{-1}$ of size $n\times n$, where $n$ is the number of vertices. 
This is a significant reduction form $\mathcal O(p^3)$  to $\mathcal O(n^3)$, when the number of edges $p$ is much larger than the number of vertices $n$. 

However, the complexity $\mathcal O(n^3)$  is still not desirable for large-scale graphs. In order to reduce the complexity further, we use the InADMM algorithm introduced in~\cite{yue2018implementing}. The key idea in the InADMM algorithm is to approximately solve a system of linear equations instead of evaluating the matrix inversion exactly. In particular, the InADMM uses the matrix identity~\eqref{eq:matrix-identity} to replace the $\beta$ update of the ADMM  iterations~\eqref{eq:beta-update} with   
\begin{align*}
& h^{k} = Q^{T}y+\rho(\alpha^{k} - w^{k})\\
& \eta^{k+1} = (QQ^{T}+\rho I)^{-1}Q h^{k})\\
& {\beta}^{k+1} = \frac{1}{\rho} (h^{k} - Q^{T}\eta^{k+1})
\end{align*}
and computes $\eta^{k+1}$ approximately using the conjugate gradient (CG) method~\cite{hestenes1952methods}. 

The most expansive step in the CG method is the matrix vector multiplication $(Q Q^{T} + \rho I)x$ where $x \in \mathbb R^n$. The complexity of this multiplication is of order $\mathcal O(p)$, because the weighted incidence matrix $Q$ has $2p$ nonzero elements.  
%  $\hat{H}\gamma_{k}$ ($\gamma_{k}$ is a $n\times 1$ vector), which is equal to
%	\begin{align*}
%	(\frac{1}{\rho}Q Q^{T} + I)\gamma_{k} = \frac{1}{\rho}QQ^{T} \gamma_{k} + \gamma_{k}.
%	\end{align*}
%	The weighted incidence matrix $Q$ has $2p$ nonzero elements. Hence,  the complexity of computing $Q^{T}x$ and then $Q(Q^Tx)$ is of order $\mathcal O(p)$. 
	Assuming the CG algorithm terminates  in $T_{CG}$ iterations, the complexity of the CG step of the InADMM algorithm is of order $\mathcal{O}(pT_{CG})$. It is straightforward to see that the complexity of other operations in InADMM is at most $\mathcal{O}(p)$. Table~\ref{tab:complexity} summarizes the complexity analysis of ADMM and InADMM algorithms.  
%Thus, the complexity of ADMM and InADMM per iteration is given as the table below:
\begin{table}[H]
	\centering 
	\begin{tabular}{c cc} 
		\hline \hline 
		Variables      & ADMM                   & InADMM\\ [0.5ex]
		\hline
		Cholesky       & $\mathcal{O}(n^{3}$)   &  -  \\
		$\eta$         & $\mathcal{O}(np)$      & $\mathcal{O}(pT_{CG})$ \\ 
		$\beta$        & $\mathcal{O}(p)$       & $\mathcal{O}(p)$  \\
		$\alpha$       & $\mathcal{O}(p)$       & $\mathcal{O}(p)$  \\
		$v$            & $\mathcal{O}(p)$       & $\mathcal{O}(p)$  \\[1ex] 
		\hline
	\end{tabular}
\caption{Computational complexity of ADMM and InADMM per iteration}
\label{tab:complexity}
\end{table}

%where
%\begin{align*}
%\hat{H} &= \frac{1}{\rho} QQ^{T}+I, \\
%\end{align*}

% the system of linear equations $\hat{H}\eta^{k} = \frac{1}{\rho}Qh^{k}$ in ADMM algorithm, InADMM algorithm applies conjugate gradient (CG) method~\cite{hestenes1952methods} to get an approximate solution of $\eta$ (denoted by $\tilde{\eta}$). The CG method require matrix $\hat{H}$, vector $\frac{1}{\rho}Qh^{k}$ and the tolerance $t$ of the residual $r_{k+1} = r_{k} - \alpha_{k}A\gamma_{k}$ as it's inputs. 

In the following sections, we present numerical experiments of applications of  the ADMM algorithm and InADMM algorithm for two examples. For the ADMM algorithm,  
we used the software code available at~\cite{lassocode2011} with the following choice of parameters: the augmented Lagrangian parameter $\rho = 1\times 10^{-7}$, the over-relaxation parameter $a = 1$, the tolerance of  primal norm $\tau_{p} = 10^{-5}$, the tolerance of residual norm $\tau_{d}= 10^{-4}$ and $\lambda = 1\times 10^{-8}\lambda_{max}$ where $\lambda_{max} = Q^{T}y$. For the InADMM algorithm, we used the CG method from~\cite{hestenes1952methods} with tolerance $10^{-4}$.  For more details about choosing the tolerance, which also guarantees the convergence of InADMM algorithm, see~\cite{yue2018implementing}.

\subsection{Random Graph}\label{subsec:Rand-Graph}
\begin{figure}[H]
	\centering % <-- added
	\begin{subfigure}{0.27\textwidth}
		\includegraphics[width=\textwidth]{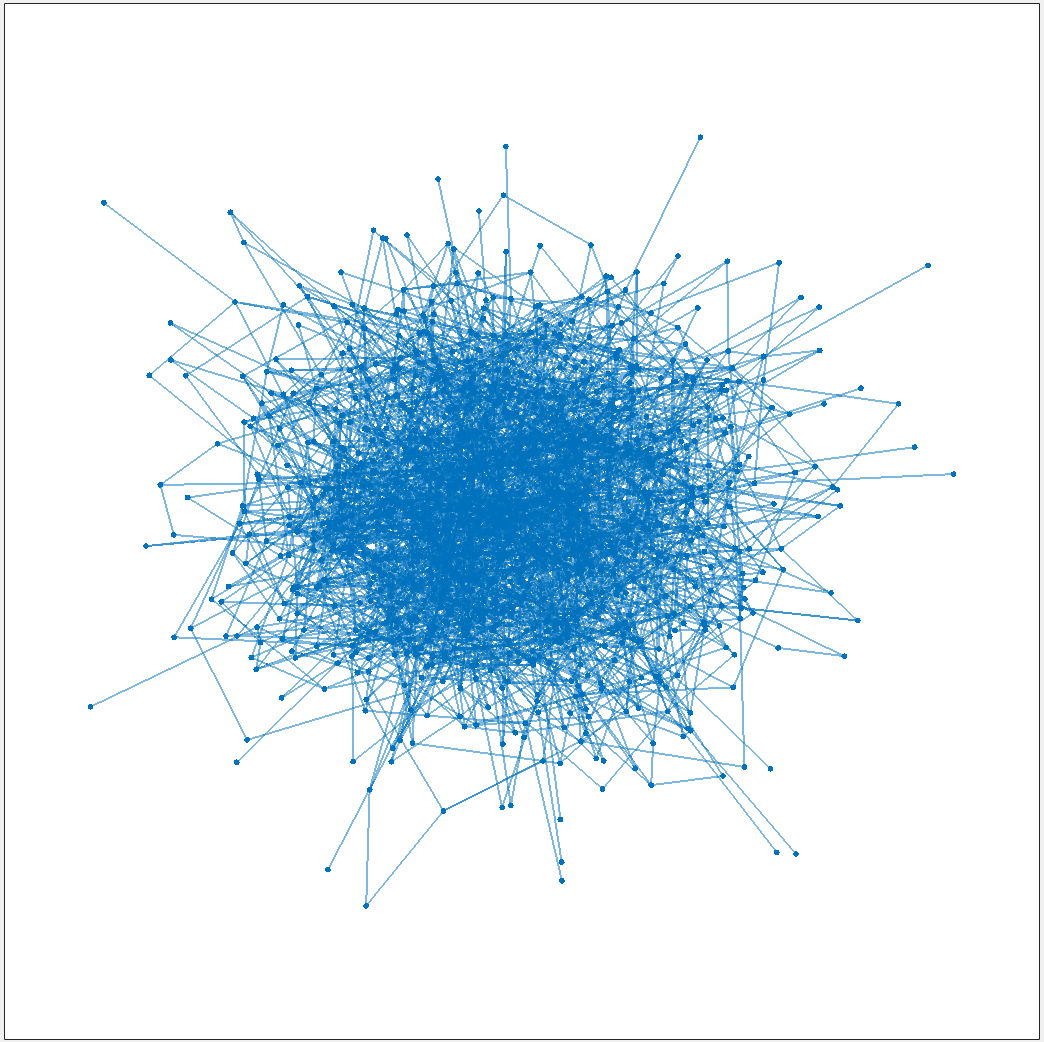}
		\caption{Random Graph}
		\label{fig:RandomGraph}
	\end{subfigure}\hfil % <-- added
	\begin{subfigure}{0.32\textwidth}
		\includegraphics[width=\linewidth]{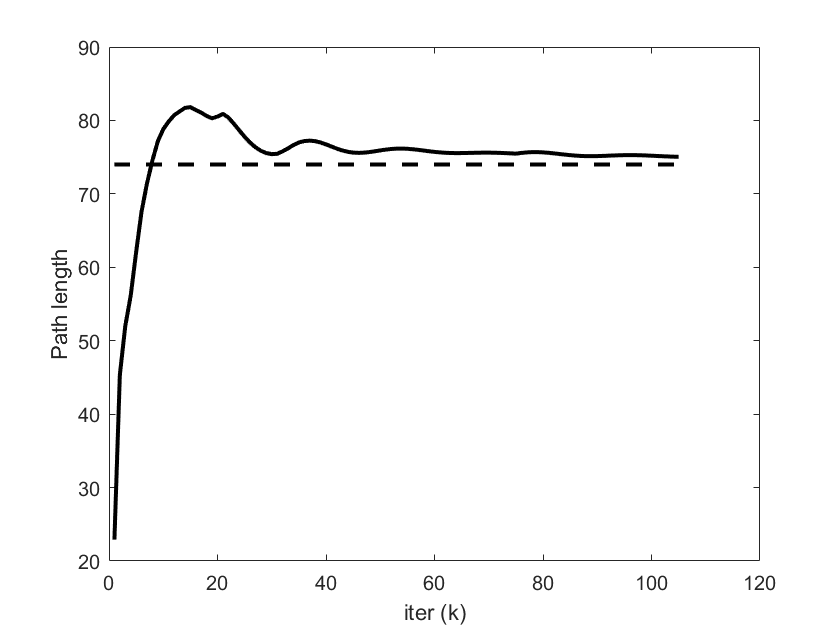}
		\caption{Estimate of the shortest path length as a function of iteration in ADMM}
		\label{fig:Convergence-lasso-RG-ADMM}
	\end{subfigure}\hfil % <-- added
    \begin{subfigure}{0.32\textwidth}
        \includegraphics[width=\linewidth]{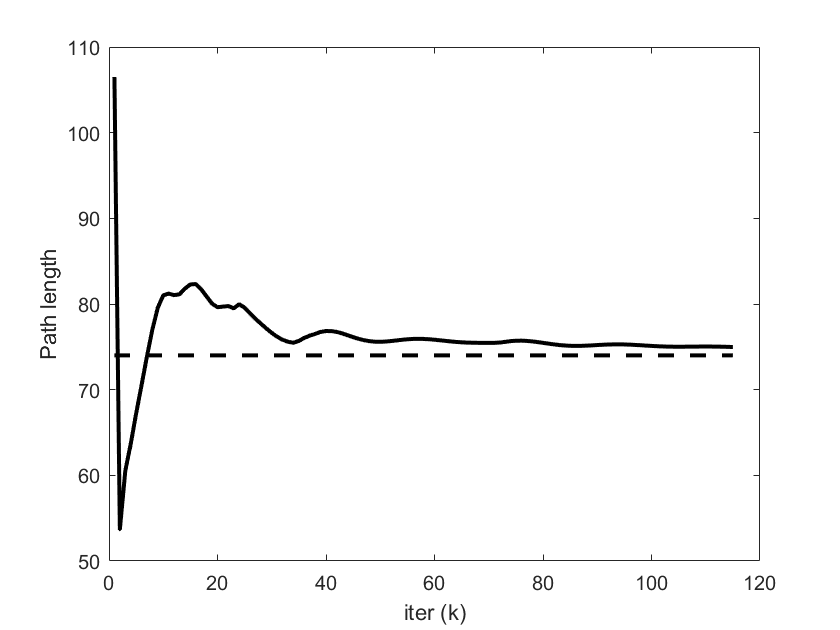}
        \caption{Estimate of the shortest path length as a function of iteration in InADMM}
        \label{fig:Convergence-lasso-RG-InADMM}
    \end{subfigure}
	\caption{Application of the ADMM and InADMM algorithm to find the length of the shortest path in a random graph with $1000$ vertices and $2688$ edges, as described in Section~\ref{subsec:Rand-Graph}}
	\label{fig:ConvergenceADMM_RandomGraph}
\end{figure}

The ADMM and InADMM algorithms are applied to find the shortest path in a random graph as depicted in Figure~\ref{fig:RandomGraph}. The random graph has $1000$ vertices. The edges are assigned randomly between two vertices with probability $2.6821\times 10^{-5}$. This yields $2688$ edges.  The weight is sampled from uniform distribution on the interval $[10,20]$. The source and target vertices are randomly picked. The result for the length of the shortest path  $\|\beta^{k}\|_1$ as a function of iterations, using the ADMM and the InADMM algorithms, is depicted in~Figure.~\ref{fig:Convergence-lasso-RG-ADMM} and~\ref{fig:Convergence-lasso-RG-InADMM}. The dashed line in the Figure represents the exact shortest path length obtained by Dijkstra's algorithm. It is observed that the lasso solution converges to the exact solution in around $50$ iterations. The running time of each iteration in InADMM in smaller than ADMM, thus the totally running time of InADMM is also smaller than InADMM (empirically $0.6$ of ADMM algorithm). 
 
\subsection{Intelligent Scissors} \label{subsec:In-Scissor}
We consider an image processing application of the shortest path problem. The application is 
 \textit{intelligent scissors} (Live-wire), which is a popular tool for image segmentation~\cite{mortensen1995intelligent}. In this application, the pixels of the image form the vertices of a graph, where each pixel is connected via an edge to its  $8$ neighbor pixels. With a suitable choice of  weights on the edges, the shortest path between two pixels is the boundary of an object~\cite[Section 3]{mortensen1995intelligent}.

\begin{figure}[H]
	\centering % <-- added
	\begin{subfigure}{0.25\textwidth}
		\includegraphics[width=\textwidth]{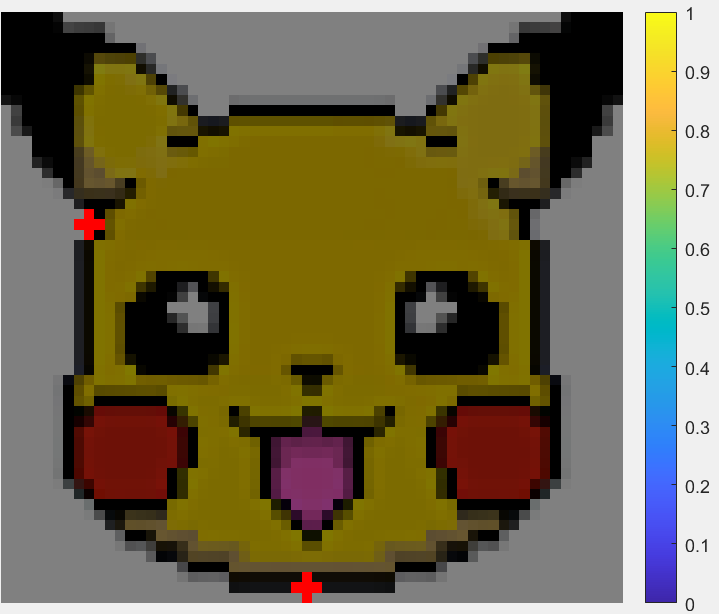}
		\caption{Original image}
		\label{fig:Org}
	\end{subfigure}\hfil % <-- added
	\begin{subfigure}{0.25\textwidth}
		\includegraphics[width=\linewidth]{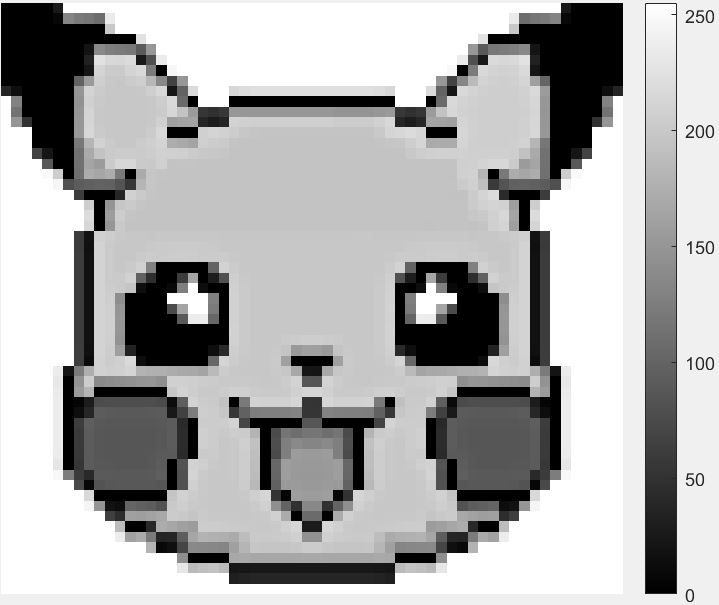}
		\caption{Grey-scale image}
		\label{fig:Grey}
	\end{subfigure}\hfil % <-- added
	\begin{subfigure}{0.25\textwidth}
		\includegraphics[width=\linewidth]{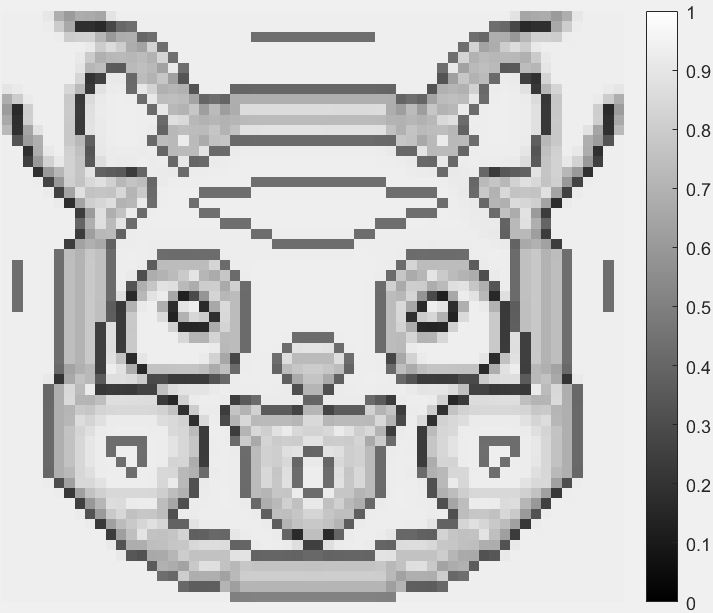}
		\caption{Edge weights}
		\label{fig:local-cost}
	\end{subfigure}

    \medskip
    \centering
	\begin{subfigure}{0.24\textwidth}
		\includegraphics[width=\textwidth]{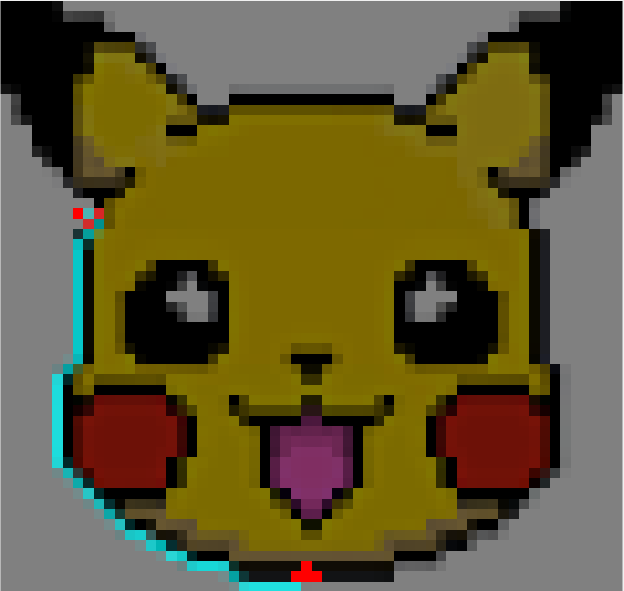}
		\caption{Dijkstra's solution path}
		\label{fig:Dijk-pika}
	\end{subfigure}\hfil % <-- added
	\begin{subfigure}{0.24\textwidth}
		\includegraphics[width=\linewidth]{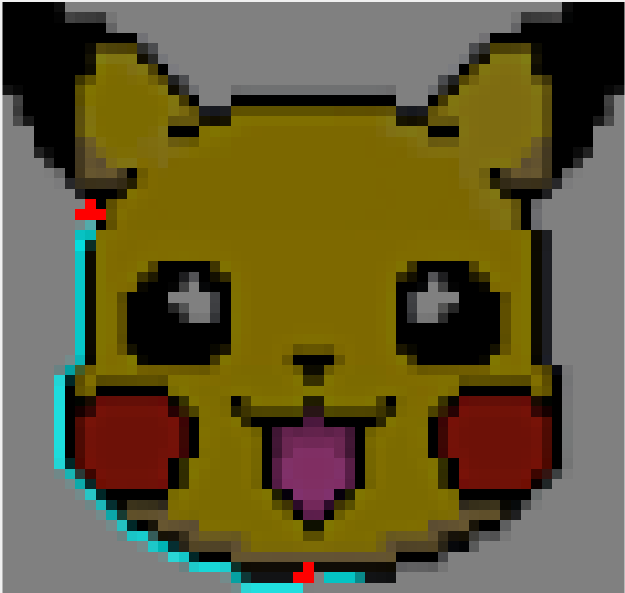}
		\caption{ADMM solution path}
		\label{fig:AMDD-lasso-pika}
	\end{subfigure}\hfil % <-- added
	\begin{subfigure}{0.24\columnwidth}
		\includegraphics[width=\linewidth]{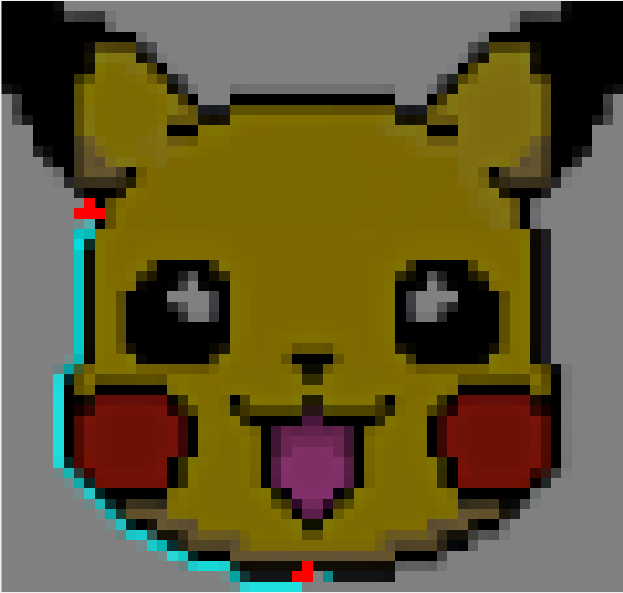}
		\caption{InADMM solution path}
		\label{fig:InADMM-lasso-pika}
	\end{subfigure}
	\caption{Application of the ADMM algorithm to find the shortest path in the \textit{Intelligent Scissor} problem as described in Section~\ref{subsec:In-Scissor}}
\end{figure}

 We apply the intelligent scissors to the Pikachu image shown in~\ref{fig:Org}. The gray-scale of the image and the structure of the edge weights between the pixels are depicted in Figure~\ref{fig:Grey} and~\ref{fig:local-cost} respectively. The picture contains $3420$ pixels (vertices) which yield $13331$ edges. The objective is to distinguish a clear boundary between the Pikachu icon and background. The objective is formulated as finding the shortest path  from pixel $(16,6)$ to pixel $(56,30)$ (top left corner to the bottom of Pikachu) as shown in Figure~\ref{fig:Org}.  
 
\begin{figure}[H]
	\centering % <-- added
	\begin{subfigure}{0.4\textwidth}
		\includegraphics[width=\textwidth]{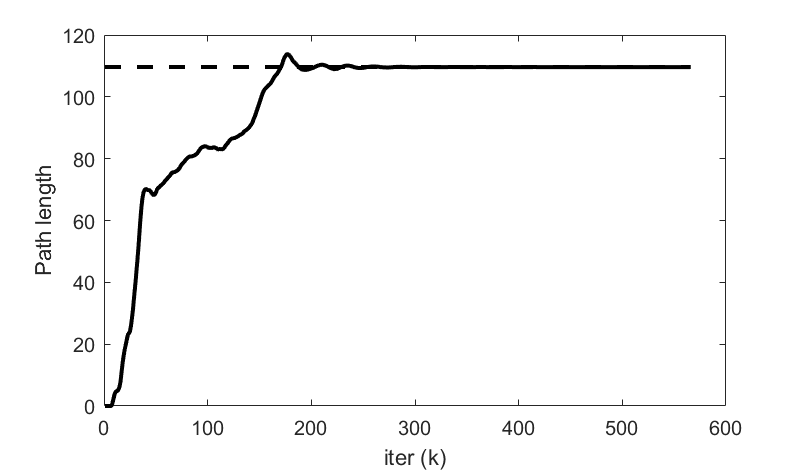}
		\caption{ADMM convergence plot}
		\label{fig:ConvergenceADMM-pika-ADMM}
	\end{subfigure}\hfil % <-- added
	\begin{subfigure}{0.4\textwidth}
		\includegraphics[width=\linewidth]{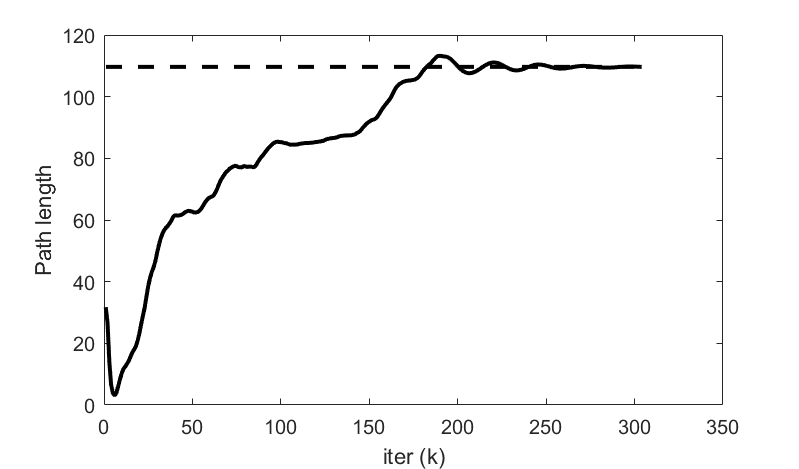}
		\caption{InADMM convergence plot}
		\label{fig:ConvergenceADMM-pika-InADMM}
	\end{subfigure}
	\caption{Estimate of the shortest path length as a function of iteration in ADMM and InADMM}
	\label{fig:Convergence-InADMM}
\end{figure}

The ADMM algorithm is simulated for this task with the same parameters as before. The resulting shortest path and the convergence of the length of the path are shown in \ref{fig:AMDD-lasso-pika} and~\ref{fig:ConvergenceADMM-pika-ADMM} respectively. For comparison, the exact shortest path obtained by the Dijkstra's algorithm is depicted in Figure~\ref{fig:Dijk-pika}. It is observed  that the ADMM algorithm provides an approximate path very similar to the exact path in around $300$ iterations. 
% The shortest path between the starting point (seed pixel) and the targeting point (pixel) is then the boundary of two different according to the property of the local cost matrix. 

As for the InADMM algorithm, the tolerance  in the CG method  is set to $10^{-7}$. 
% It shows that by choosing the $t$ which is closed to it upper bound, the value of $T$ can be efficiently reduced. 
 The resulting shortest path and the convergence of the length of the path are shown in \ref{fig:InADMM-lasso-pika} and~\ref{fig:ConvergenceADMM-pika-InADMM} respectively

\section{Conclusion}
In this paper, we studied the Lasso formulation of the shortest path problem. We showed that the lasso path solution is equivalent to the shortest path trees that appear in the Dijkstra's algorithm. And we  proposed to apply the ADMM algorithm to estimate the shortest path length for large graphs. 
Careful analysis of the computational complexity and the distributed implementation of the ADMM, for this particular objective, is subject of future work. 

\appendix
%\balance

\subsection{Derivation of \eqref{eq:crossingtime}} \label{apdx:t-cross-derivation}
For simplicity, we drop the iteration subscript $k$ in our derivations. 
$D_\calA$ is the incidence matrix formed by the edges in the active set. The graph formed by $\calA$ consist of two disjoint trees $\Ts$, $\Tt$, and set of isolated vertices $\Omega$. We decompose the rows of matrix $D_\calA$ into rows corresponding to these three subsets, and express $D_\calA$ according to
\begin{align*}
D_{\calA} = \begin{bmatrix}
D_{\calA^{(s)}} & {0} \\
0 &  0 \\
0 &  D_{\calA^{(t)}}
\end{bmatrix}
\end{align*}
where $D_{\calA^{(s)}}$ and $D_{\At}$ are the incidence matrix for the tree $\Ts$ and $\Tt$ respectively, and $0$ represents the zero matrix of appropriate dimensions. Then, 
\begin{align*}
D_\calA^+ = \begin{bmatrix}
D_{\calA^{(s)}}^+ & 0 & 0 \\
0 & 0 & D_{\calA^{(t)}}^+  
\end{bmatrix}
\end{align*}
We use this expression and Lemma~\ref{lem:incidence} to compute $a$ and $b$. By definition~\eqref{eq:a-b}
\begin{align*}
a &=(Q_{\calA}^{T}Q_{\mathcal{A}})^{+}Q_{\mathcal{A}}^{T}y 
= Q_A^+y=W_\calA D_\calA^+y \\
&=\begin{bmatrix}
-\frac{1}{|\Ts|}W_{\As} P^{(s)} \mathbbm{1}_s \\
\frac{1}{|\Tt|}W_{\At} P^{(t)} \mathbbm{1}_t
\end{bmatrix} 
%b_j &= (Q_\calA^TQ_\calA)^+ s = W_{\mathcal{A}}D_{\mathcal{A}}^{+} (D_{\mathcal{A}}^{T})^{+}W_{\mathcal{A}}s
%& = W_{\mathcal{A}}    
%\begin{bmatrix}
%-\frac{1}{n}P \mathbbm{1} & PJ
%\end{bmatrix}
%\begin{bmatrix}
%(-\frac{1}{n} P \mathbbm{1})^{T} \\
%(PJ)^{T}
%\end{bmatrix} 
%W_{\mathcal{A}}s\\
%%& = W_{\mathcal{A}} \Big(\frac{1}{n^2}P \mathbbm{1} \mathbbm{1}^{T} P +
%%P M_{1,1} M_{1,1}^{T} P\Big)W_{\mathcal{A}} s\\
%%& = W_{\mathcal{A}}\Big(\frac{1}{n^2}P \mathbbm{1} \mathbbm{1}^{T} P +
%%P (I-\3frac{1}{n})(I-\frac{1}{n})^{T} P\Big)W_{\mathcal{A}} s\\
%& =-W_{\mathcal{A}}\Big(P \calL -\frac{1}{n}P \mathbbm{1}  \mathbbm{1}^T\calL)
\end{align*}
where $P^{(s)}$ and $P^{(t)}$ are the path matrix for tree $\Ts$ and $\Tt$ respectively, and $\mathbbm{1}_t$ and $\mathbbm{1}_s$ are all one vectors of size $|\Ts|$ and $|\Tt|$ respectively. For $b$,
\begin{align*}
%a_j &=(Q_{\calA}^{T}Q_{\mathcal{A}})^{+}Q_{\mathcal{A}}^{T}y 
%= Q_A^+y=W_\calA D_\calA^+y \\
%&=\begin{bmatrix}
%-\frac{1}{|\Ts|}W_{\As} P^{(s)} \mathbbm{1}^{(s)} \\
%\frac{1}{|\Tt|}W_{\At} P^{(t)} \mathbbm{1}^{(t)}
%\end{bmatrix} 
b&= (Q_\calA^TQ_\calA)^+ s = W_{\mathcal{A}}D_{\mathcal{A}}^{+} (D_{\mathcal{A}}^{T})^{+}W_{\mathcal{A}}s
%\\& = W_{\mathcal{A}}    
%\begin{bmatrix}
%-\frac{1}{n}P \mathbbm{1} & PJ
%\end{bmatrix}
%\begin{bmatrix}
%(-\frac{1}{n} P \mathbbm{1})^{T} \\
%(PJ)^{T}
%\end{bmatrix} 
%W_{\mathcal{A}}s\\
%& = W_{\mathcal{A}} \Big(\frac{1}{n^2}P \mathbbm{1} \mathbbm{1}^{T} P +
%P M_{1,1} M_{1,1}^{T} P\Big)W_{\mathcal{A}} s\\
%& = W_{\mathcal{A}}\Big(\frac{1}{n^2}P \mathbbm{1} \mathbbm{1}^{T} P +
%P (I-\3frac{1}{n})(I-\frac{1}{n})^{T} P\Big)W_{\mathcal{A}} s\\
\\& = \begin{bmatrix}-W_{\As}\Big(P^{(s)} \calL^{(s)} -\frac{1}{|\Ts|}P^{(s)} \mathbbm{1}_s \mathbbm{1}_s^T\calL^{(s)}\Big)\\
W_{\At}\Big(P^{(t)} \calL^{(t)} -\frac{1}{|\Tt|}P^{(t)} \mathbbm{1}_t  \mathbbm{1}_t^T\calL^{(t)}\Big)
\end{bmatrix}
\end{align*}
%Hence, we decompose the rows of $D_\calA$ into these three subsets of corresponding to vertices in the tree rooted at $s$, isolated vertices, and the tree rooted at $t$.   As a result, some vertices are isolated. 
where $\calL^{(s)} \triangleq -(P^{(s)})^T\calW_{\calA^{(s)}} s$ is a vector of size $|\Ts|$ corresponding to vertices in the tree $\Ts$. The component of $\calL^{(s)}$,  corresponding to vertex $v \in \Ts$, is equal to $l_v^{(s)}$, i.e. the length of the path from  $v$ to the root $s$. The vector $\calL^{(t)} \triangleq (P^{(t)})^T\calW_{\calA^{(t)}} s$ has similar interpretation, but for vertices of tree $\Tt$. 

Putting the results for $a$ and $b$ together, the ratio $a_j/b_j$ for $e_j \in \calA^{(s)}$ is
\begin{align*} 
\frac{a_j}{b_j}&= \frac{\frac{1}{|\Ts|}w_j|R_j|}{w_j(\sum_{v\in R_j}l_v^{(s)} - \frac{|R_j|}{|\Ts|} \sum_{v\in \Ts} l_v^{(s)})}
%\\
%&= 
%\frac{-\frac{1}{n}W_{\mathcal{A}} P \mathbbm{1}}{W_{\mathcal{A}}(P P^{T} -\frac{1}{n}P J P^{T})W_{\mathcal{A}} s} \\
%&= \frac{\frac{1}{n}W_{\mathcal{A}} P \mathbbm{1}}{W_{\mathcal{A}} (P\mathcal{L} -\frac{1}{n}P J \mathcal{L})}
\end{align*}
where $R_j$ is the set of non-zero components of the $j$-th row of $P^{(s)}$.  This concludes~\eqref{eq:crossingtime} for $e_j \in \calA^{(s)}$. The derivation for $e_j \in \calA^{(t)}$ is similar.

\subsection{Derivation of~\eqref{eq:tjoin}} \label{apdx:t-join-derivation}
By definition of joining time~\eqref{eq:join}
\begin{align}\label{eq:t-join-apdx}
t^{\text{join}}_j = \frac{\frac{1}{w_j}D_j^T(Q_\calA a - y)}{\frac{1}{w_j}D_j^T(Q_\calA b) \pm 1}
\end{align}
Next, we obtain expressions for the terms in parentheses. For the term in the numerator
\begin{align*}
Q_\calA a - y &= D_\calA D_\calA^+ y - y \\
& =  \begin{bmatrix}
D_{\calA^{(s)}} & {0} \\
0 &  0 \\
0 &  D_{\calA^{(t)}}
\end{bmatrix}
\begin{bmatrix}
D_{\calA^{(s)}}^+ & {0} & 0 \\
0 & 0 &  D_{\calA^{(t)}}^+
\end{bmatrix}y - y \\
&= \begin{bmatrix}
-\frac{1}{|\Ts|}\mathbbm{1}_s \\
0 \\
+\frac{1}{|\Tt|}\mathbbm{1}_t
\end{bmatrix}
\end{align*}
where we used $D D^+=I - \frac{1}{\mathbbm{1}^T\mathbbm{1}}\mathbbm{1}\mathbbm{1}^T$ for each incidence matrix $D=D_{\As}$ and $D=D_{\At}$. And for the term in the denominator
\begin{align*}
Q_\calA b &= D_\calA D_\calA^+ (D_\calA^T)^+ W_\calA s = 
(D_\calA^+)^T W_\calA s \\
&=\begin{bmatrix}- \calL^{(s)} + \frac{1}{|\Ts|} \mathbbm{1}_s \mathbbm{1}_s^T\calL^{(s)}\\ 0 \\
\calL^{(t)} -\frac{1}{|\Tt|} \mathbbm{1}_t  \mathbbm{1}_t^T\calL^{(t)}
\end{bmatrix}
\end{align*}
Using these results in~\eqref{eq:t-join-apdx} and evaluating the expression for $e_j$ for each case in~\eqref{eq:tjoin} concludes~\eqref{eq:tjoin}. 

\subsection{Proof of lemma 3.2} \label{apdx:proof-lemma-3.2}
The proof is based on \cite[lemma 2]{tibshirani2013lasso}. The active set $\calA$ is always unique. 
In Section~\ref{sec:relation}, we showed that the active set form two disjoint trees. Hence,
% From the connection between Dijkstra's algorithm and LARS algorithm, we know that:\\
%For $\lambda_{k} > 0$, the trees with root $v_{s}$ and $v_{t}$ are not connect,
$Q_{\mathcal{A}} = \begin{bmatrix}
Q_{\mathcal{A}^{(s)}} &Q_{\mathcal{A}^{(t)}}
\end{bmatrix}$
where $Q_{\mathcal{A}^{(s)}} = D_{\mathcal{A}_{}^{(s)}} W_{\mathcal{A}_{}^{(s)}}^{-1}$ 
and $Q_{\mathcal{A}_{}^{(t)}} = D_{\mathcal{A}_{}^{(t)}} W_{\mathcal{A}_{}^{(t)}}^{-1}$.
$D_{\mathcal{A}_{k}^{(s),(t)}}$ are incidence matrices of two trees and $W_{\mathcal{A}_{k}^{(s).(t)}}$ are two diagonal matrices  with positive elements. The kernel of incidence matrix of a tree is empty because there is no cycle. Hence, the rank is equal to the number of columns. As a result, the rank of $Q_\calA$ is equal to $|\calA|$. Then, according to~\cite[lemma 2]{tibshirani2013lasso} the Lasso solution is unique.   
%equal to the   Thus, every column of $Q_{\mathcal{A}_{k}^{(s)}}$ are independent to each other and every column of $Q_{\mathcal{A}_{k}^{(s)}}$ are also independent to each other. Since $\mathcal{A}_{k}^{(s)}$ and $\mathcal{A}_{k}^{(t)}$ are two disconnect trees, it can be concluded that every column of $Q_{\mathcal{A}_{k}}$ are independent to each other, or equivalently rank($Q_{\mathcal{A}_{k}}$)$= |\mathcal{A}_{k}|$. Then from \cite[lemma 2]{tibshirani2013lasso}, the lasso solution $\beta_{\calA_k}$ is unique.\\
%For $\lambda_{k} = 0$, the two trees $T^{(s)}$ and $T^{(t)}$ are connect and form one tree. By the Assumption A1 that the shortest path is unique, the lasso solution $\beta_{\calA_k}$ is the shortest path and $Q_{\mathcal{A}_{k}^{(s)}}\beta_{\calA_k} = y^{(s,t)}$ satisfied.\\
%In conclusion, for $\lambda_{k} \geq 0$, the lasso solution $\beta_{\calA_k}$ is always unique.
\section*{Acknowledgments}
Partial funding was provided by NSF under grants 1807664, 1839441, AFOSR under grant FA9550-20-1-0029.
\bibliographystyle{IEEEtran}
\bibliography{ref}    
\end{document}